\documentclass{amsart}
\usepackage{amssymb}
\usepackage{amscd}
\usepackage{verbatim}
\usepackage{epsfig}
\begin{document}
\newcommand\Mand{\ \text{and}\ }
\newcommand\Mwith{\ \text{with}\ }
\newcommand\Mfor{\ \text{for}\ }
\newcommand\Mst{\ \text{such that}\ }
\newcommand\Mor{\ \text{or}\ }
\newcommand\Mif{\ \text{if}\ }
\newcommand\Miff{\ \text{iff}\ }
\newcommand\Mthen{\ \text{then}\ }
\newcommand\nin{\notin}
\newcommand\identity{\operatorname{id}}
\newcommand\Id{\operatorname{Id}}
\newcommand\Real{\mathbb{R}}
\newcommand\pos{\Real^+}
\newcommand\Rnp{\Real\setminus\{0\}}
\newcommand\nzero{\setminus\{0\}}
\newcommand\Cx{\mathbb{C}}
\newcommand\Cxp{\Cx^+}
\newcommand\Cxm{\Cx^-}
\newcommand\Nat{\mathbb{N}}
\newcommand\halfNat{{\frac{1}{2}}\mathbb{N}}
\newcommand\intgr{\mathbb{Z}}
\newcommand\im{\operatorname{Im}}
\newcommand\re{\operatorname{Re}}
\newcommand\sign{\operatorname{sign}}
\newcommand\codim{\operatorname{codim}}
\newcommand\End{\operatorname{End}}
\newcommand\Ker{\operatorname{Ker}}
\newcommand\Hom{\operatorname{Hom}}
\newcommand\tr{\operatorname{tr}}
\newcommand\Tr{\operatorname{Tr}}
\newcommand\ideal{{\mathcal I}}
\newcommand\Span{\operatorname{span}}
\newcommand\image{\operatorname{image}}
\newcommand\Range{\operatorname{Ran}}
\newcommand\Graph{\operatorname{graph}}
\newcommand\slim{\operatornamewithlimits{s-lim}}
\newcommand\spp{\operatorname{sp}}
\newcommand\sll{\operatorname{sl}}
\newcommand\sol{\operatorname{so}}
\newcommand\SL{\operatorname{SL}}
\newcommand\SO{\operatorname{SO}}
\newcommand\On{\operatorname{O}}
\newcommand\pa{\partial}
\newcommand\Rn{\Real^n}
\newcommand\Rm{\Real^m}
\newcommand\RN{\Real^N}
\newcommand\RtN{\Real^{2N}}
\newcommand\RM{\Real^M}
\newcommand\sphere{\mathbb{S}}
\newcommand\Sn{\sphere^{n-1}}
\newcommand\Sm{\sphere^{m-1}}
\newcommand\Snp{\sphere^n_+}
\newcommand\Smp{\sphere^m_+}
\newcommand\SN{\sphere^{N-1}}
\newcommand\SNp{\sphere^N_+}
\newcommand\circlep{\sphere^1_+}
\newcommand\Phom{P_{h}}
\newcommand\Shom{S_{h}}
\newcommand\distance{\operatorname{dist}}
\newcommand\cl{\operatorname{cl}}
\newcommand\interior{\operatorname{int}}
\newcommand\Fa{\operatorname{Fa}}
\newcommand\ff{\operatorname{ff}}
\newcommand\mf{\operatorname{mf}}
\newcommand\cf{\operatorname{cf}}
\newcommand\scf{\operatorname{sf}}
\newcommand\lf{\operatorname{lf}}
\newcommand\rf{\operatorname{rf}}
\newcommand\indfam{{\mathcal K}}
\newcommand\fraka{{\mathfrak a}}
\newcommand\calA{{\mathcal A}}
\newcommand\calB{{\mathcal B}}
\newcommand\calH{{\mathcal H}}
\newcommand\calR{{\mathcal R}}
\newcommand\calO{{\mathcal O}}
\newcommand\calJ{{\mathcal J}}
\newcommand\calM{{\mathcal M}}
\newcommand\calN{{\mathcal N}}
\newcommand\calX{{\mathcal X}}
\newcommand\calF{{\mathcal F}}
\newcommand\calG{{\mathcal G}}
\newcommand\calT{{\mathcal T}}
\newcommand\calP{{\mathcal P}}
\newcommand\calPt{\widetilde{\mathcal P}}
\newcommand\calS{{\mathcal S}}
\newcommand\calSt{\widetilde{\mathcal S}}
\newcommand\calC{{\mathcal C}}
\newcommand\calCt{{\tilde {\mathcal C}}}
\newcommand\calCL{{\mathcal C}_{\text L}}
\newcommand\calCR{{\mathcal C}_{\text R}}
\newcommand\Cinf{{\mathcal C}^{\infty}}
\newcommand\dist{{\mathcal C}^{-\infty}}
\newcommand\dCinf{\dot\Cinf}
\newcommand\ddist{\dot\dist}
\newcommand\Cj{{\mathcal C}^j}
\newcommand\Linf{L^{\infty}}
\newcommand\phg{{\text{phg}}}
\newcommand\comp{{\text{comp}}}
\newcommand\bcon{{\mathcal A}}
\newcommand\bconc{{\mathcal A}_{\text{phg}}}
\newcommand\Sch{{\mathcal S}}
\newcommand\temp{\Sch^{\prime}}
\newcommand\Diff{\operatorname{Diff}}
\newcommand\Diffb{\operatorname{Diff}_{\text{b}}}
\newcommand\Diffc{\operatorname{Diff}_{\text{c}}}
\newcommand\Diffsc{\operatorname{Diff}_{\text{sc}}}
\newcommand\DiffI{\operatorname{Diff}_{\text{I}}}
\newcommand\DiffIq{\operatorname{Diff}_{\text{I},q}}
\newcommand\sing{\text{sing}}
\newcommand\reg{\text{reg}}
\newcommand\supp{\operatorname{supp}}
\newcommand\ssupp{\operatorname{sing\ supp}}
\newcommand\csupp{\operatorname{cone\ supp}}
\newcommand\esupp{\operatorname{ess\ supp}}
\newcommand\Fr{{\mathcal F}}
\newcommand\Frinv{\Fr^{-1}}
\newcommand\bop{{\mathcal B}}
\newcommand\spec{\operatorname{spec}}
\newcommand\pspec{\spec_{pp}}
\newcommand\cspec{\spec_{c}}
\newcommand\FIO{{\mathcal I}}
\newcommand\SP{\operatorname{RC}}
\newcommand\RC{\operatorname{RC}}
\newcommand\Symc{S_c}
\newcommand\Symca{S_c^{\alpha}}
\newcommand\Symczero{S_c^{0,...,0}}
\newcommand\sci{{}^{\text{sc}}}
\newcommand\sct{\sci T^*}
\newcommand\scT{\sci T}
\newcommand\scdt{\sci \dot T^*}
\newcommand\dS{\dot S^*}
\newcommand\dT{\dot T^*}
\newcommand\dSreg{\dot\Sigma_{\text reg}}
\newcommand\scct{\sci\bar{T}^*}
\newcommand\Csc{C_{\text{sc}}}
\newcommand\SNpscd{(\SNp)^2_{\text{sc}}}
\newcommand\scdiag{\Delta_{\text{sc}}}
\newcommand\projscl{\pi^L_{\text{sc}}}
\newcommand\projscr{\pi^R_{\text{sc}}}
\newcommand\scHL{\sci H^{2,0}_{|\zeta|^2-\lambda^2}}
\newcommand\scHrg{\sci H^{2,0}_{\sqrt{g}}}
\newcommand\Hsc{H_{\text{sc}}}
\newcommand\WF{\operatorname{WF}}
\newcommand\WFp{\operatorname{WF^{\prime}}}
\newcommand\WFsc{\operatorname{WF}_{\text{sc}}}
\newcommand\WFscp{\operatorname{WF_{sc}^{\prime}}}
\newcommand\WFC{\operatorname{WF}_C}
\newcommand\WFCi{\operatorname{WF}_{C_i}}
\newcommand\elliptic{\operatorname{ell}}
\newcommand\Psop{\operatorname{\Psi}}
\newcommand\Psiscrs{\operatorname{\Psi_{sc}^{-2,\infty}}}
\newcommand\Psiscr{\operatorname{\Psi_{sc}^{-2,0}}}
\newcommand\Psiscrm{\operatorname{\Psi_{sc}^{0,2}}}
\newcommand\PsiscHam{\operatorname{\Psi_{sc}^{2,0}}}
\newcommand\Psisci{\operatorname{\Psi_{sc}^{*,*}}}
\newcommand\Psiscid{\operatorname{\Psi_{sc}^{0,0}}}
\newcommand\Psiscis{\operatorname{\Psi_{sc}^{0,\infty}}}
\newcommand\Psiscsi{\operatorname{\Psi_{sc}^{-\infty,0}}}
\newcommand\Psiscs{\operatorname{\Psi_{sc}^{-\infty,\infty}}}
\newcommand\Psiscalg{\operatorname{\Psi_{sc}^{\infty,-\infty}}}
\newcommand\nullHam{{\mathcal N}}
\newcommand\charD{\Sigma_{\Delta-\lambda^2}}
\newcommand\charLap{\Sigma_{\Delta-\lambda}}
\newcommand\Snl{\Sn_{\lambda}}
\newcommand\SNl{\SN_{\lambda}}
\newcommand\gammat{\tilde\gamma}
\newcommand\gammasc{\gamma}
\newcommand\Tau{\mathcal{T}}
\newcommand\taut{\tilde\tau}
\newcommand\taub{\bar\tau}
\newcommand\Nout{N^+_{\lambda}}
\newcommand\Nin{N^-_{\lambda}}
\newcommand\Nio{N^{\pm}_{\lambda}}
\newcommand\El{E_{\lambda}}
\newcommand\Elt{\tilde E_{\lambda}}
\newcommand\Eil{E^i_{\lambda}}
\newcommand\Ejl{E^j_{\lambda}}
\newcommand\Eajl{E^{\alpha_j}_{\lambda}}
\newcommand\Eilt{\tilde E^i_{\lambda}}
\newcommand\Np{N^+}
\newcommand\Nm{N^-}
\newcommand\Npm{N^{\pm}}
\newcommand\Fin{F^-(\lambda)}
\newcommand\Fini{F^-_i(\lambda)}
\newcommand\Fout{F^+(\lambda)}
\newcommand\Fouti{F^+_i(\lambda)}
\newcommand\Foutj{F^+_j(\lambda)}
\newcommand\Rout{R^+_{\lambda}}
\newcommand\Routl{R^+_{\lambda^2}}
\newcommand\Routsgnl{R^{\sign\lambda}_{\lambda^2}}
\newcommand\Rin{R^-_{\lambda}}
\newcommand\Rinl{R^-_{\lambda^2}}
\newcommand\Rinsgnl{R^{-\sign\lambda}_{\lambda^2}}
\newcommand\Rio{R^{\pm}_{\lambda}}
\newcommand\Riol{R^{\pm}_{\lambda^2}}
\newcommand\Roi{R^{\mp}_{\lambda}}
\newcommand\Roil{R^{\mp}_{\lambda^2}}
\newcommand\Riob{R^{\pm}}
\newcommand\Roib{R^{\mp}}
\newcommand\Tio{T^{\pm}}
\newcommand\Tiob{T^{\pm}_{\ff}}
\newcommand\Toi{T^{\mp}}
\newcommand\Toib{T^{\mp}_{\ff}}
\newcommand\TIiob{T_I^{\pm}}
\newcommand\Rinb{R^-}
\newcommand\Rinbsgnl{R^{-\sign\lambda}}
\newcommand\Tin{T^-}
\newcommand\Tinb{T^-_{\ff}}
\newcommand\TIinb{T^-_I}
\newcommand\Routb{R^+}
\newcommand\Routbsgnl{R^{\sign\lambda}}
\newcommand\Tout{T^+}
\newcommand\Toutb{T^+_{\ff}}
\newcommand\TIoutb{T^+_I}
\newcommand\Rlkf{(|\xib|^2-(\lambda-i0)^2)^{-1}}
\newcommand\Rlk{\rho_0(\lambda)}
\newcommand\Rmlk{\rho_0(-\lambda)}
\newcommand\Rpmlk{\rho_0(\pm\lambda)}
\newcommand\Rlka{\rho_1(\lambda)}
\newcommand\Rlkb{\rho_2(\lambda)}
\newcommand\Rilk{\rho_i(\lambda)}
\newcommand\reduced{\natural}
\newcommand\Rlf{R_0(\lambda)}
\newcommand\Rla{R_1(\lambda)}
\newcommand\Rlb{R_2(\lambda)}
\newcommand\Ril{R_i(\lambda)}
\newcommand\Rlj{R_j(\lambda)}
\newcommand\Rlft{R_0(\lambda)}
\newcommand\Rflambda{R_0^{\reduced}(\sigma)}
\newcommand\RV{R^{\reduced}_V}
\newcommand\Rfsigma{R_0^{\reduced}(\sigma)}
\newcommand\Rfsigmah{R_0^{\reduced}(\sigma^{1/2})}
\newcommand\Rfzero{R_0^{\reduced}(0)}
\newcommand\RlV{R^{\reduced}_V(\sigma)}
\newcommand\RlVi{R^{\reduced}_{V_i}(\sigma)}
\newcommand\RlVt{R_V(\lambda)}
\newcommand\RlVtL{{R}_V^L(\lambda)}
\newcommand\RlVtR{{R}_V^R(\lambda)}
\newcommand\RlVit{{R}_{V_i}(\lambda)}
\newcommand\RlVta{{R}_V^{(1)}(\lambda)}
\newcommand\RlVtk{{R}_V^{(k)}(\lambda)}
\newcommand\RlVatV{{R}_{V_{\alpha}}(\lambda)V_{\alpha}}
\newcommand\RlVatVa{{R}_{V_{\alpha_1}}(\lambda)V_{\alpha_1}}
\newcommand\RlVatVb{{R}_{V_{\alpha_2}}(\lambda)V_{\alpha_2}}
\newcommand\RlVatVk{{R}_{V_{\alpha_k}}(\lambda)V_{\alpha_k}}
\newcommand\RlVatVkk{{R}_{V_{\alpha_{k+1}}}(\lambda)V_{\alpha_{k+1}}}
\newcommand\RlVaptV{{R}_{V_{\alpha'}}(\lambda)V_{\alpha'}}
\newcommand\RlVapptV{{R}_{V_{\alpha''}}(\lambda)V_{\alpha''}}
\newcommand\RlVajtV{{R}_{V_{\alpha_j}}(\lambda)V_{\alpha_j}}
\newcommand\RlVaktV{{R}_{V_{\alpha_k}}(\lambda)V_{\alpha_k}}
\newcommand\RlVakktV{{R}_{V_{\alpha_{k+1}}}(\lambda)V_{\alpha_{k+1}}}
\newcommand\Tl{T(\lambda)}
\newcommand\Tlt{\tilde\Tl}
\newcommand\Tltp{\tilde T'(\lambda)}
\newcommand\Tltpp{\tilde T''(\lambda)}
\newcommand\Tli{T_i(\lambda)}
\newcommand\Tlit{\tilde\Tli}
\newcommand\Tlip{T_i'(\lambda)}
\newcommand\Tlipp{T_i''(\lambda)}
\newcommand\Tlj{T_j(\lambda)}
\newcommand\Tla{T_{\alpha}(\lambda)}
\newcommand\Tlaa{T_{\alpha_1}(\lambda)}
\newcommand\Tlab{T_{\alpha_2}(\lambda)}
\newcommand\Tlak{T_{\alpha_k}(\lambda)}
\newcommand\Tlakt{\tilde\Tlak}
\newcommand\Tlaj{T_{\alpha_j}(\lambda)}
\newcommand\Tlajj{T_{\alpha_{j+1}}(\lambda)}
\newcommand\Tlajp{T_{\alpha_j}'(\lambda)}
\newcommand\Tlajpt{\tilde\Tlajp}
\newcommand\Tlajt{\tilde\Tlaj}
\newcommand\Tlakk{T_{\alpha_{k+1}}(\lambda)}
\newcommand\Tlakkp{T_{\alpha_{k+1}}'(\lambda)}
\newcommand\Tlap{T_{\alpha'}(\lambda)}
\newcommand\Tlapt{\tilde\Tlap}
\newcommand\Tlapp{T_{\alpha''}(\lambda)}
\newcommand\Tkl{T^{(k)}(\lambda)}
\newcommand\Tcl{T^{\flat}(\lambda)}
\newcommand\Fl{F(\lambda)}
\newcommand\BlVt{\tilde B_V(\lambda)}
\newcommand\KBlVt{K_{\BlVt}}
\newcommand\BlVaat{B_{V_{\alpha_1}}(\lambda)}
\newcommand\BV{B_V}
\newcommand\Bone{B_1}
\newcommand\Btwo{B_2}
\newcommand\Bthree{B_3}
\newcommand\Banyj{B_j}
\newcommand\PlV{P_V(\lambda)}
\newcommand\PlVc{P_V^{\flat}(\lambda)}
\newcommand\Pl{P_0(\lambda)}
\newcommand\SVl{S_V(\lambda)}
\newcommand\Sjr{S_j^{\reduced}}
\newcommand\Rkp{{\mathcal R}^k_+}
\newcommand\Rkm{{\mathcal R}^k_-}
\newcommand\Rkpm{{\mathcal R}^k_{\pm}}
\newcommand\Phys{{\mathcal P}}
\newcommand\Pc{\overline{\mathcal P}}
\newcommand\pip{\pi^{\perp}}
\newcommand\pipa{\pi_\partial}
\newcommand\gammapa{\gamma_\partial}
\newcommand\pipah{\hat\pi_\partial}
\newcommand\pit{\tilde\pi}
\newcommand\xit{\tilde\xi}
\newcommand\zetat{\tilde\zeta}
\newcommand\etat{\tilde\eta}
\newcommand\sigmat{\tilde\sigma}
\newcommand\sigmahat{\hat\sigma}
\newcommand\thetat{\tilde\theta}
\newcommand\psit{\tilde\psi}
\newcommand\phit{\tilde\phi}
\newcommand\chit{\tilde\chi}
\newcommand\rhot{\tilde\rho}
\newcommand\xib{\bar\xi}
\newcommand\zetab{\bar\zeta}
\newcommand\thetab{\bar\theta}
\newcommand\etab{\bar\eta}
\newcommand\iotal{\iota_{\lambda}}
\newcommand\rhoat{\rhot_{\alpha_1}}
\newcommand\Lambdat{\tilde\Lambda}
\newcommand\Lambdati{\tilde\Lambda^{\text{in}}}
\newcommand\Lambdato{\tilde\Lambda^{\text{out}}}
\newcommand\Lambdatp{\tilde\Lambda^{\text{prop}}}
\newcommand\Lambdai{\Lambda^{\text{in}}}
\newcommand\Lambdao{\Lambda^{\text{out}}}
\newcommand\poles{\Lambda'}
\newcommand\rpoles{\Lambda_p}
\newcommand\thresholds{\Lambda}
\newcommand\Vt{\tilde V}
\newcommand\It{\tilde I}
\newcommand\half{{\frac{1}{2}}}
\newcommand\sigmah{\sigma^{1/2}}
\newcommand\bX{\partial X}
\newcommand\bXb{\partial \Xb}
\newcommand\Deltabt{\tilde\Delta_0}
\newcommand\strip{\Omega_T}
\newcommand\Kf{K^{\flat}}
\newcommand\Gs{G^{\sharp}}
\newcommand\Gt{\tilde G}
\newcommand\Osc{\sci\Omega}
\newcommand\OSc{{}^\Scl\Omega}
\newcommand\Osch{\sci\Omega^{\half}}
\newcommand\Oscmh{\sci\Omega^{-\half}}
\newcommand\Isc{I_{sc}}
\newcommand\os{{\text{os}}}
\newcommand\Qzl{Q^0_{-\lambda}}
\newcommand\Lie{{\mathcal L}}
\newcommand\bl{{\text b}}
\newcommand\scl{{\text{sc}}}
\newcommand\sccl{{\text{scc}}}
\newcommand\Scl{{\text{Sc}}}
\newcommand\ScLl{{\text{Sc,L}}}
\newcommand\ScRl{{\text{Sc,R}}}
\newcommand\Sccl{{\text{Scc}}}
\newcommand\sus{{\text{sus}}}
\newcommand\ssl{{\text{ss}}}
\newcommand\XXb{X^2_\bl}
\newcommand\XXbt{\Xt^2_\bl}
\newcommand\XXsc{X^2_\scl}
\newcommand\XXsct{\Xt^2_\scl}
\newcommand\XXSc{X^2_\Scl}
\newcommand\XXSct{\Xt^2_\Scl}
\newcommand\XXScL{X^2_\ScLl}
\newcommand\XXScR{X^2_\ScRl}
\newcommand\MMsc{M^2_\scl}
\newcommand\Deltab{\Delta_\bl}
\newcommand\Deltasc{\Delta_\scl}
\newcommand\DeltaSc{\Delta_\Scl}
\newcommand\DeltaScL{\Delta_\ScLl}
\newcommand\DeltaScR{\Delta_\ScRl}
\newcommand\prs{\sigma}
\newcommand\Nsc{N_\scl}
\newcommand\Nscp{N_{\scl,p}}
\newcommand\Nff{N_{\ff}}
\newcommand\Nffz{N_{\ff,0}}
\newcommand\Nffzp{N_{\ff,0,p}}
\newcommand\Nffl{N_{\ff,l}}
\newcommand\Nffml{N_{\ff,-l}}
\newcommand\Nmf{N_{\mf}}
\newcommand\Nmfz{N_{\mf,0}}
\newcommand\Nmfl{N_{\mf,l}}
\newcommand\Nmfml{N_{\mf,-l}}
\newcommand\ffb{\operatorname{bf}}
\newcommand\Ffb{\operatorname{bf'}}
\newcommand\ffsc{\operatorname{sf}}
\newcommand\ffSc{\operatorname{sf_C}}
\newcommand\Ffsc{\operatorname{sf'}}
\newcommand\rff{\rho_{\ff}}
\newcommand\rmf{\rho_{\mf}}
\newcommand\rffb{\rho_{\ffb}}
\newcommand\rffsc{\rho_{\ffsc}}
\newcommand\rFfsc{\rho_{\Ffsc}}
\newcommand\rffSc{\rho_{\ffSc}}
\newcommand\rinf{\rho_{\infty}}
\newcommand\CL{C_L}
\newcommand\CR{C_R}
\newcommand\betab{\beta_\bl}
\newcommand\betasc{\beta_\scl}
\newcommand\betaSc{\beta_\Scl}
\newcommand\BetaSc{\bar\beta_\Scl}
\newcommand\betaScL{\beta_\ScLl}
\newcommand\betaScR{\beta_\ScRl}
\newcommand\ScT{{}^\Scl T^*}
\newcommand\SccT{{}^\Scl \bar T^*}
\newcommand\ScS{{}^\Scl S^*}
\newcommand\Tb{{}^\bl T}
\newcommand\Tsc{{}^\scl T}
\newcommand\TSc{{}^\Scl T}
\newcommand\CSc{C_\Scl}
\newcommand\Lambdasc{{}^\scl\Lambda}
\newcommand\XXXb{X^3_\bl}
\newcommand\XXXsc{X^3_\scl}
\newcommand\XXXSc{X^3_\Scl}
\newcommand\XXXScO{X^3_{\Scl,O}}
\newcommand\XXXScF{X^3_{\Scl,F}}
\newcommand\XXXScS{X^3_{\Scl,S}}
\newcommand\XXXScC{X^3_{\Scl,C}}
\newcommand\KDsc{\operatorname{KD^{\half}_\scl}}
\newcommand\KDSc{\operatorname{KD^{\half}_\Scl}}
\newcommand\KDScEF{\operatorname{KD^{E,F}_\Scl}}
\newcommand\Oh{\operatorname{\Omega^{\half}}}
\newcommand\WFSc{\WF_\Scl}
\newcommand\WFtSc{\WF_{\text 3sc}}
\newcommand\WFScmf{\WF_{\Scl,\mf}}
\newcommand\WFScff{\WF_{\Scl,\ff}}
\newcommand\WFScs{\WF_{\Scl,\prs}}
\newcommand\WFScp{\WF'_\Scl}
\newcommand\WFScmfp{\WF'_{\Scl,\mf}}
\newcommand\WFScffp{\WF'_{\Scl,\ff}}
\newcommand\WFScsp{\WF'_{\Scl,\prs}}
\newcommand\Diffscc{\Diff_\sccl}
\newcommand\DiffSc{\Diff_\Scl}
\newcommand\DiffScc{\Diff_\Sccl}
\newcommand\DiffscI{\Diff_{\scl,\text{I}}}
\newcommand\VscI{\Vf_{\scl,\text{I}}}
\newcommand\DiffsV{\operatorname{Diff}_{\sus(V)}}
\newcommand\DiffsVsc{\operatorname{Diff}_{\sus(V),\scl}}
\newcommand\DiffsVCsc{\operatorname{Diff}_{\sus(V)-C,\scl}}   
\newcommand\Psisc{\Psop_\scl}
\newcommand\Psiscc{\Psop_\sccl}
\newcommand\Psiss{\Psop_\ssl}
\newcommand\Psisch{\Psop_{\scl,h}}
\newcommand\Psiscch{\Psop_{\sccl,h}}
\newcommand\PsiSc{\Psop_\Scl}
\newcommand\PsiScph{\Psop_{\Scl,\phi}}
\newcommand\PsiScra{\Psop_{\Scl,\rho^\sharp_a}}
\newcommand\PsiScc{\Psop_\Sccl}
\newcommand\PsiSccml{\Psop^{m,l}_\Sccl}
\newcommand\PsiScxx{\Psop^{*,*}_\Scl}
\newcommand\PsiScml{\Psop^{m,l}_\Scl}
\newcommand\PsiScmz{\Psop^{m,0}_\Scl}
\newcommand\PsiScmmz{\Psop^{-m,0}_\Scl}
\newcommand\PsiSckz{\Psop^{k,0}_\Scl}
\newcommand\PsiScmmml{\Psop^{-m,-l}_\Scl}
\newcommand\Psiscmkk{\Psop^{-k,k}_\scl}
\newcommand\Psiscmmmkk{\Psop^{-m-k,k}_\scl}
\newcommand\Psiscmoo{\Psop^{-1,1}_\scl}
\newcommand\Psiscmz{\Psop^{m,0}_\scl}
\newcommand\Psiscmmz{\Psop^{-m,0}_\scl}
\newcommand\PsiSckmkl{\Psop^{km,kl}_\Scl}
\newcommand\PsiScmplp{\Psop^{m',l'}_\Scl}
\newcommand\PsiScmmpllp{\Psop^{m+m',l+l'}_\Scl}
\newcommand\Psiscml{\Psop^{m,l}_\scl}
\newcommand\PsiScid{\Psop^{0,0}_\Scl}
\newcommand\PsiSczo{\Psop^{0,1}_\Scl}
\newcommand\PsiScmii{\Psop^{-\infty,\infty}_\Scl}
\newcommand\PsiScmiz{\Psop^{-\infty,0}_\Scl}
\newcommand\PsiScmoo{\Psop^{-1,1}_\Scl}
\newcommand\PsisCid{\Psop^{0,0}_{\scl-C}}
\newcommand\PsisC{\Psop_{\scl-C}}
\newcommand\Psiinf{\Psop_{\infty}}
\newcommand\Psiinfid{\Psop_{\infty}^0}
\newcommand\PsiFinf{\Psop_{\infty-\Fr}}
\newcommand\PsisVscml{\Psop^{m,l}_{\sus(V),\scl}}
\newcommand\PsisVsc{\Psop_{\sus(V),\scl}}
\newcommand\PsisVpsc{\Psop_{\sus(V_p),\scl}}
\newcommand\PsisVCSc{\Psop_{\sus(V)-C,\scl}}
\newcommand\SFinf{S_{\infty-\Fr}}
\newcommand\YsVC{Y^2_{\sus(V)-C,\scl}}
\newcommand\ffYsc{\ffsc_{\sus(V)}}
\newcommand\SXC{S(X;C)}
\newcommand\Ios{I_{\text{os}}}
\newcommand\pbL{\pi^2_{\bl,{\text L}}}
\newcommand\pbR{\pi^2_{\bl,{\text R}}}
\newcommand\pscL{\pi^2_{\scl,{\text L}}}
\newcommand\pscR{\pi^2_{\scl,{\text R}}}
\newcommand\PbO{\pi^3_{\bl,{\text O}}}
\newcommand\PscO{\pi^3_{\scl,{\text O}}}
\newcommand\PScO{\pi^3_{\Scl,{\text O}}}
\newcommand\PScF{\pi^3_{\Scl,{\text F}}}
\newcommand\PScC{\pi^3_{\Scl,{\text C}}}
\newcommand\PScS{\pi^3_{\Scl,{\text S}}}
\newcommand\pScL{\pi^2_{\Scl,{\text L}}}
\newcommand\pScR{\pi^2_{\Scl,{\text R}}}
\newcommand\CLF{\CL^F}
\newcommand\CLO{\CL^O}
\newcommand\CLS{\CL^S}
\newcommand\CLC{\CL^C}
\newcommand\DeltaYb{\Delta_{\bl,Y}}
\newcommand\DeltaYsc{\Delta_{\sus-\scl}}
\newcommand\diag{\operatorname{diag}}
\newcommand\Vf{{\mathcal V}}
\newcommand\Vb{{\mathcal V}_{\bl}}
\newcommand\Vsc{{\mathcal V}_{\scl}}
\newcommand\VSc{{\mathcal V}_{\Scl}}
\newcommand\VfI{\Vf_{\text{I}}}
\newcommand\VfIq{\Vf_{\text{I},q}}
\newcommand\scH{{}^\scl H}
\newcommand\scHg{\scH_g}
\newcommand\Hss{H_\ssl}
\newcommand\xh{\hat x}
\newcommand\yh{\hat y}
\newcommand\sh{\hat s}
\newcommand\rh{\hat r}
\newcommand\Yh{\hat Y}
\newcommand\Zh{\hat Z}
\newcommand\Yb{\bar Y}
\newcommand\hb{\bar h}
\newcommand\xih{\hat\xi}
\newcommand\etah{\hat\eta}
\newcommand\muh{\hat\mu}
\newcommand\mub{\bar\mu}
\newcommand\nub{\bar\nu}
\newcommand\mubh{\widehat{\bar\mu}}
\newcommand\yb{\bar y}
\newcommand\zb{\bar z}
\newcommand\ub{\bar u}
\newcommand\Qb{\bar Q}
\newcommand\Wbp{{\bar W}^\perp}
\newcommand\Wp{W^\perp}
\newcommand\Kt{\tilde K}
\newcommand\Wt{\tilde W}
\newcommand\Ut{\tilde U}
\newcommand\yt{\tilde y}
\newcommand\ut{\tilde u}
\newcommand\vt{\tilde v}
\newcommand\ft{\tilde f}
\newcommand\htil{\tilde h}
\newcommand\St{\tilde S}
\newcommand\Pt{\tilde P}
\newcommand\Rt{\tilde R}
\newcommand\qt{\tilde q}
\newcommand\Qt{\tilde Q}
\newcommand\Xb{\bar X}
\newcommand\lambdat{\tilde\lambda}
\newcommand\betat{\tilde\beta}
\newcommand\Phit{\tilde\Phi}
\newcommand\epst{\tilde\epsilon}
\newcommand\ep{\epsilon}
\newcommand\bt{\tilde b}
\newcommand\Xt{\tilde X}
\newcommand\Mt{\tilde M}
\newcommand\At{\tilde A}
\newcommand\Et{\tilde E}
\newcommand\Ht{\tilde H}
\newcommand\at{\tilde a}
\newcommand\Ct{\tilde C}
\newcommand\pih{\hat\pi}
\newcommand\Rh{\hat R}
\newcommand\Ah{\hat A}
\newcommand\Bh{\hat B}
\newcommand\Ch{\hat C}
\newcommand\Gh{\hat G}
\newcommand\Hh{\hat H}
\newcommand\Qh{\hat Q}
\newcommand\Ph{\hat P}
\newcommand\Nh{\hat N}
\newcommand\Sh{\hat S}
\newcommand\Gcal{{\mathcal G}}
\newcommand\GcalC{{\mathcal G}_C}
\newcommand\Jcal{{\mathcal J}}
\newcommand\JcalC{{\mathcal J}_C}
\newcommand\evpr{\lambda_1}
\newcommand\evth{\lambda_0}
\setcounter{secnumdepth}{3}
\newtheorem{lemma}{Lemma}[section]
\newtheorem{prop}[lemma]{Proposition}
\newtheorem{thm}[lemma]{Theorem}
\newtheorem{cor}[lemma]{Corollary}
\newtheorem{result}[lemma]{Result}
\newtheorem*{thm*}{Theorem}
\newtheorem*{prop*}{Proposition}
\newtheorem*{conj*}{Conjecture}
\numberwithin{equation}{section}
\theoremstyle{remark}
\newtheorem{rem}[lemma]{Remark}
\theoremstyle{definition}
\newtheorem{Def}[lemma]{Definition}
\newtheorem*{Def*}{Definition}
\def\signature#1#2{\par\noindent#1\dotfill\null\\*
{\raggedleft #2\par}}

\renewcommand{\theenumi}{\roman{enumi}}
\renewcommand{\labelenumi}{(\theenumi)}

\title[Low energy inverse problems]
{Low energy inverse problems in three-body scattering}
\author[Gunther Uhlmann and Andras Vasy]{Gunther Uhlmann and Andr\'as Vasy}
\date{October 15, 2001}
\address{Department of Mathematics, University of Washington, 
Seattle, WA}
\email{gunther@math.washington.edu}
\address{Department of Mathematics, Massachusetts Institute of Technology,
Cambridge MA 02139}
\email{andras@math.mit.edu}
\thanks{G.\ U.\ is partially supported by NSF grant \#DMS-00-70488 and
a John Simon Guggenheim fellowship.
A.\ V.\ is partially supported by NSF grant \#DMS-99-70607. Both authors
are grateful for the hospitality of the Erwin Schr\"odinger Institute
in Vienna and the Mathematical Sciences Research Institute in Berkeley, CA}

\maketitle

\section{Introduction}
Scattering theory is the analysis of the motion of several interacting
particles.
Inverse problems in scattering theory seek the answer to the question:
can one determine the interactions between particles by a scattering
experiment, so, for instance, does the scattering matrix, or a part of the
scattering matrix, determine the interaction? The answer, and the difficulty,
greatly depends on the part of the scattering data one wishes to use.
Before explaining the various settings, we remark that the inverse
problems are highly non-linear since the scattering data do not
depend linearly on the interactions. Hence, one of the usual methods
in the field is to transfer the problem to an asymptotically linear
one, which is then easier to analyze.

The simplest setting is the study of high energy asymptotics of the
scattering matrix, for then the potentials behave like small perturbations
of the Laplacian. In addition, the leading term in the asymptotics depends
linearly on them. This problem
was studied, under various assumptions, by Enss and Weder
\cite{Enss-Weder:Geometrical, Enss-Weder:Inverse},
Novikov \cite{Novikov:N-body} and Wang \cite{Wang:High}.

Here we are interested in finite energy problems, i.e.\ where
the scattering data are known either only at a fixed energy, or in a fixed
bounded interval of energies. Thus, the problems are not immediately
equivalent to a linear perturbation problem. The flavor of the
problem greatly depends on the part of the scattering matrix one
wishes to use.

In some situations a principal symbol
calculation for an S-matrix allows one to use the 2-body inverse
results. An example of this is free-to-free scattering: as shown in
\cite{Vasy:Structure}, in three-body scattering
the singularities of the free-to-free
S-matrix at energy $\lambda>0$ determine the S-matrices in all proper
subsystems at all energies in $(0,\lambda)$, which then determine
the pair interactions by two-body results. More precisely, the principal
symbol of the part of the free-to-free S-matrix
corresponding to a single collision is essentially given by the
subsystem S-matrix at the energies in $(0,\lambda)$. Slightly more
involved arguments using the results of \cite{Vasy:Bound-States} are
expected to work
in the many-body setting to show that the free-to-free S-matrix determines
all pair interactions.

However, one may wish to study inverse problems where the parts of
the S-matrix that are known do not have any singularities, so the
previous method cannot be applied. An example is two-cluster to any other
cluster scattering. Indeed, Skibsted \cite{Skibsted:Smoothness} has shown that
the corresponding S-matrices have smooth kernels
apart from the diagonal singularity of the 2-cluster to same 2-cluster
S-matrix. The latter agrees with the diagonal singularity of the kernel
of the identity operator if the potentials are Schwartz, hence is of
no help for the inverse problem; this is also the case with
the two-body problem with Schwartz potentials. In certain ways,
these are the most realistic problems, for in a scattering experiment
one typically shoots a particle at a nucleus, atom,
molecule, or other such composite
`cluster', which may break up as a result of the collision. One then
measures the outcome of the experiment -- this is exactly the information
contained in the two-cluster to other-cluster S-matrices.

We study two-cluster to same two-cluster scattering for three-body
Hamiltonians with real-valued potentials
under the assumption that all unknown interactions are short-range and small.
We show that the S-matrices $S_{\alpha\alpha'}(\lambda)$
in an energy interval $I_\lambda$ below the break-up energy
determine the Fourier transform of the effective
interaction in a ball, whose radius is determined by the energy of the bound
state under consideration. More precisely we prove the following theorem,
whose statement uses some notation that we describe in detail in the
next section.

\begin{thm}\label{thm:main}
Suppose that $\alpha$ is a channel in a 2-cluster $a$, $\dim X_a\geq 2$,
and $\mu>\dim X_a$, and $V_a\in S^{-s}(X^a)$
for some $s>0$.
There exists a constant $\delta>0$ such that the following statement holds.

Suppose that $\sup|\langle w^b\rangle^\mu V_b|<\delta$ for all $b\neq a$.
Suppose also that $I\subset(\ep_\alpha,0)$ is
a non-empty open interval, and
let $R=2\sqrt{\sup I-\ep_\alpha}$. Then $S_{\alpha'\alpha''}(\lambda)$,
$\lambda\in I$, given for all bound states $\alpha'$, $\alpha''$ with energy
$\ep_{\alpha'},\ep_{\alpha''}<0$, determines the Fourier transform
$\hat V_\alpha$ of the effective interaction,
\begin{equation}\label{eq:pair-101-intro}
V_\alpha=\int_{X^a}
I_a|\psi_{\alpha}|^2\,dw^a,
\end{equation}
in the ball $B_0(R)$ of radius $R$
centered at the origin in $X_a$.
\end{thm}

If we only want to determine $\hat V_\alpha$ in a smaller ball, we
need even less information. There is a variety of statements one
can make using different information; we only make the following one.

\begin{thm}\label{thm:main2}
Suppose that $a$ is a 2-cluster, $\dim X_a\geq 2$,
$V_a\in S^{-s}(X^a)$ for some $s>0$,
$\alpha$ is the ground state of $H^a$
with $\ep_\alpha<0$.
Let $\ep'>\ep_\alpha$ be the next eigenvalue of $H^a$, or $0$ if this
does not exist. Let $\mu>\dim X_a$.
There exists a constant $\delta>0$ such that the following statement holds.

Suppose that $\sup|\langle w^b\rangle^\mu V_b|<\delta$, for all $b\neq a$.
Suppose also that $I\subset(\ep_\alpha,0)$ is
a non-empty open interval, and
let $R=2\sqrt{\min(\sup I,\ep')-\ep_\alpha}$. Then $S_{\alpha\alpha}(\lambda)$,
$\lambda\in I$, determines the Fourier transform
$\hat V_\alpha$ of the effective interaction,
\begin{equation}\label{eq:pair-101-intro-p}
V_\alpha=\int_{X^a}
I_a|\psi_{\alpha}|^2\,dw^a,
\end{equation}
in the ball $B_0(R)$ of radius $R$
centered at the origin in $X_a$.
\end{thm}

Since we are working below the break-up energy, heuristically one
expects that the composite particle may be regarded as a single
particle, and two-body methods may be applied. This turns out to be
false, at least when taken literally. Indeed, many two-body methods,
one of which we describe below, rely
on allowing large complex momenta for the particles, which in turn
permits the break-up of a cluster. Hence, one of the themes of this
paper is the extent to which composite particles may be regarded as a
single unit below break-up energies for the purposes of inverse
problems in scattering theory.

Our strategy is similar to how one approaches low energy inverse problems
in two-body scattering, which we now briefly recall.
Then the kernel of the (relative) S-matrix is also smooth for Schwartz
potentials, once the kernel of $\Id$ is subtracted, and is conormal
to the diagonal for symbolic potentials.
Faddeev \cite{Faddeev:Increasing, Faddeev:Factorization, Faddeev:Inverse}
started the study of exponential solutions, i.e.\ solutions
of $(H-\lambda)u=0$ of the form $u=u_\rho=
e^{i\rho\cdot w}(1+v_\rho(w))$, $v_\rho$
`small', $\rho$ not necessarily real, and $\rho\cdot\rho=\lambda$.
Even if $\lambda$ is fixed, by allowing $\rho$ to be complex, one
can take $\rho\to\infty$, so that $v_\rho\to 0$ in an appropriate sense.
Provided that one can relate the pairing
\begin{equation}\label{eq:pairing-int}
\int u_\rho V e^{-i\rho'\cdot w}\,dw,
\end{equation}
taking the limit $\rho\to\infty$ (and $\rho'\to\infty$)
becomes an analogue of the high energy
limit, with the leading term linear in the potential. In other words,
the high energy asymptotics is replaced by high complex momentum
asymptotics,
as pioneered by Calder\'on,
see \cite{Calderon:Inverse, Sylvester-Uhlmann:Global, RBMGeo}.

While the S-matrix is not
analytic in the energy $\lambda$ unless other assumptions are made,
for a very large class of potentials (including Schwartz potentials)
$u_\rho$ is meromorphic (indeed, analytic if $V$ is small)
in the complex one-dimensional space (i.e.\ line) spanned by
$\im\rho$, provided this line is fixed. In other words, $u_\rho$ is
analytic in $z$ (in $\im z\neq 0$), where we write $\rho=z\nu+\rho_\perp$,
$\nu,\rho_\perp$ real, $\nu\cdot\rho_\perp=0$. Moreover, from $\im z>0$,
$u_\rho$ extends continuously to $\im z=0$ if $V$ is small, and to a large
subset of the real line in $z$ otherwise. This can be exploited
in problems where the S-matrix is known in an interval, as in the
work of Novikov \cite{Novikov-Khenkin:D-bar}, Weder \cite{Weder:Global}
and Isozaki~\cite{Isozaki:Multi-dimensional}.
Indeed, one shows first that the S-matrix in an energy interval determines
the pairing \eqref{eq:pairing-int} in a corresponding interval in $z$,
then uses that the boundary values of a meromorphic function determine
the function, finally lets $z\to\infty$ and uses the high momentum limit
to determine the Fourier transform of $V$.

In the three-body setting, there are similar exponential solutions
corresponding to a bound state $\psi_\alpha$ of a subsystem $H^a$,
so $(H^a-\ep_\alpha)\psi_\alpha=0$. Namely, one considers solutions of
$(H-\lambda)u_\rho=0$ of the form
$u_\rho=e^{i\rho\cdot w_a}(\psi_\alpha(w^a)+v_\rho(w))$ where $\rho$ is in the
complexification $\Cx(X_a)$ of $X_a$,
$\rho\cdot\rho=\lambda-\ep_\alpha$, $\lambda\in\Cx$.
In fact, this construction works in great generality,
though the structure of $u_\rho$ changes with
$\rho$. In this paper we keep $|\rho_\perp|<\sqrt{-\ep_\alpha}$, in
which case $u_\rho$ can be constructed by perturbation theory. In
particular, it is easy to see that $u_\rho$ depends analytically on $z$.
Here perturbation theory is understood loosely, for even if the
unknown interactions $V_b$ are small, they are {\em not} a compact
perturbation of $H_a=\Delta+V_a$, for they do {\em not} decay at
infinity.
In particular, if $V_b$ becomes large,
the structure of $u_\rho$ changes drastically, and its
analyticity in $z$ is far from clear. Even for small $V_b$, if we take
$\rho_\perp$ large, the cluster will be allowed to break up, creating
a major difficulty for fixed energy inverse problems.

On the other hand, the connection to the S-matrices
is less immediate than in the two-body setting. In general, one expects
that all parts of the S-matrix need to be known at a certain energy to
determine the pairing $\int u_\rho I_a \overline{\psi_{\alpha}(w^a)}
e^{-i\rho'\cdot w_a}\,dw$. This can be seen explicitly from the
statement of our main theorem, when it is applicable:
$S_{\alpha'\alpha''+}(\lambda)$ play a role in the statement for all
$\alpha'$, $\alpha''$. In fact, analogously to an observation of Novikov
\cite{Novikov:Determination}, by reducing $|\rho_\perp|$ further (than
$|\rho_\perp|<\sqrt{-\ep_\alpha}$ mentioned above),
some of the two-cluster to two-cluster S-matrices can be eliminated,
as was done in the second theorem. 
However, the restriction on $\rho_\perp$
implies restrictions on the frequencies at which $\hat V_\alpha$ can
be recovered.

The structure of this paper is the following. After recalling the
usual many-body notation, we construct the exponential eigenfunctions,
and we study their limit as $\rho$ becomes real. We use this to relate
the corresponding pairing to the S-matrix. Finally, we apply this to
the study of the inverse problem by taking $\rho\to\infty$, and prove
Theorems~\ref{thm:main} and \ref{thm:main2}.

The authors are grateful to Rafe Mazzeo, Richard Melrose, Roman Novikov
and Maciej Zworski for helpful discussions.

\section{Notation and preliminaries}
Below the notation is that of \cite{Vasy:Bound-States}, which is
to say it is the standard many-body notation as in
\cite{Derezinski-Gerard:Scattering}. First, $X=X_0=\Rn$ is the total
configuration space, equipped with the standard Euclidean metric $g$.
The collision planes, $X_a$, $a\in I$, $I$ finite,
are linear subspaces of $X_0$, and $X^a$ is the orthocomplement
of $X_a$ in $X_0$. We assume that $\{X_a:\ a\in I\}$ is closed under
intersections, includes $X_0$ and $X_1=\{0\}$.
We write $w=(w_a,w^a)$ for coordinates on $X=X_a\oplus X^a$,
and identify $X_a^*$ with $X_a$ via the metric $g$.

We write $\sphere_a=C_a$ for the unit sphere in $X_a$ (with respect to
the metric inherited from $X_0$). Geometrically it is better to consider
$C_a$ as `the sphere at infinity', but for the sake of simplicity (and
to conform with the usual many-body conventions) we adopt the unit sphere
point of view. We also let
$C_{a,\sing}=\cap_{C_b\subsetneq C_a} C_b$
be the singular, $C_{a,\reg}=C_a\setminus C_{a,\sing}$ the regular
part of $C_a$. Recall also that a two-cluster $a$, denoted by $\#a=2$,
is a (non-trivial) cluster such that
$C_b\subset C_a$ implies $b=a$
(or $b=1$ or $a=1$ provided that $C_1=\emptyset$ is included in the collection
of $C_c$'s). Thus, two-clusters are the most singular
clusters, and
in particular if $a$ is a 2-cluster, then $C_{a,\sing}=\emptyset$.
The collection of collision plans corresponds to a three-body geometry
if every cluster, except $0$ and $1$, is a two-cluster.

Concerning the analytic aspects, we write $L^2_p(X_a)$ for the weighted
$L^2$ space $L^2(X_a,\langle w_a\rangle^{2p}\,dw_a)$ on $X_a$. We also
write $H^s_p(X_a)=H^{s,p}(X_a)$ for the Sobolev space corresponding to
this weight.

We let $H^a$ be the subsystem Hamiltonian on $X^a$,
i.e.\
\begin{equation*}H^a=\Delta_{X^a}+\sum_{X^b\subset X^a} V_b,\end{equation*}
and $I_a$ is
the intercluster interaction $I_a=V-\sum_{X^b\subset X^a} V_b$.
The unreduced subsystem Hamiltonian acts on functions on the
whole space $\Rn$; it is
\begin{equation*}H_a=\Delta_{X_a}+H^a.\end{equation*}
In addition, $R^a$, resp.\ $R_a$, denote the resolvent of the reduced,
resp.\ unreduced, Hamiltonian
of the subsystem $a$, i.e.\ $R^a(\sigma)=(H^a-\sigma)^{-1}$,
$R_a(\sigma)=(H_a-\sigma)^{-1}$
for $\sigma\nin\Real$. We write $\Lambda$ for the set of thresholds of $H$,
which is defined inductively over the proper subsystems by
\begin{equation*}
\Lambda_a=\cup_{X_b\supsetneq X_a} \Lambda'_b,
\quad \Lambda'_b=\Lambda_b\cup\pspec(H^b),
\end{equation*}
and we usually denote the spectral parameter by $\lambda$. In particular,
for a three-body Hamiltonian $H$, if $a$ is a 2-cluster then
$\Lambda_a=\{0\}$, $\Lambda'_a=\{0\}\cup\pspec(H^a)$, and $\Lambda_1=\Lambda
=\cup_{\#a=2}\Lambda'_a$, $\Lambda'=\Lambda_1'=\Lambda\cup\pspec(H)$.

Let $\psi_\alpha$ denote the normalized 
$L^2$ eigenfunctions of $H^a$, and let
$\ep_\alpha$ be the bound state energy in $\psi_\alpha$:
$(H^a-\ep_\alpha)\psi_\alpha=0$.
If $\ep_\alpha$ is not an eigenvalue
of a proper subsystem of $H^a$, then $\psi_\alpha\in
e^{-\mu_\alpha|w_a|}L^2(X^a)$
for some $\mu_\alpha>0$ given by the next threshold above $\ep_\alpha$
(see \cite{FroExp}). We call such a bound state $\alpha$ a non-threshold
bound state.

Usually the Poisson operator and the scattering matrices are considered
as operators on functions on unit spheres in appropriate spaces. When
we investigate the real-frequency
behavior of the exponential solutions that we construct
in the next section, it will be convenient to consider the Poisson operators
and S-matrices as operators acting on functions on spheres of different radii.
Thus, we replace the unit sphere in $X_a$ by
the sphere $\sphere_a(\sqrt{\lambda-\ep_\alpha})$ as the parameterization
space for the Poisson operators; here
\begin{equation*}
\sphere_a(\sigma)=\{\xi_a\in X_a:\ |\xi_a|=\sigma\}.
\end{equation*}
The regular and singular parts of $\sphere_a(\sigma)$ are defined analogously
to those of $C_a$.

The thus normalized forward Poisson operator of $H_a$ in channel
$\alpha$ is given
by
\begin{equation*}\begin{split}
\calPt_{\alpha,+}(\lambda)g=c_\alpha\psi_{\alpha}(w^a)&
\int_{\sphere_a(\sqrt{\lambda-\ep_\alpha})} e^{-iw_a\cdot\omega_a}
g\,d\omega_a,\quad g\in\Cinf_c(\sphere_{a,\reg}(\sqrt{\lambda-\ep_\alpha})),\\
&c_\alpha=(\lambda-\ep_\alpha)^{\frac{m-1}4}e^{-\frac{m-1}4\pi i}
(2\pi)^{-\frac{m-1}2},\ m=\dim X_a,
\end{split}\end{equation*}
where $d\omega_a$ is
the standard measure on $\sphere_a(\sqrt{\lambda-\ep_\alpha})$ normalized
to have volume equal to that of the unit sphere.
The Poisson operator of $H$ in channel
$\alpha$ is then
\begin{equation*}\begin{split}
\calP_{\alpha,+}(\lambda)g&=\calPt_{\alpha,+}(\lambda)g
-R(\lambda+i0)((H-\lambda)
\calPt_{\alpha,+}(\lambda)g)\\
&=\calPt_{\alpha,+}(\lambda)g-R(\lambda+i0)I_a
\calPt_{\alpha,+}(\lambda)g.
\end{split}\end{equation*}
Note that if $\langle w^b\rangle^\mu V_b\in L^\infty(X^b)$, then
\begin{equation*}
\langle w_a\rangle^{p} I_a\calPt_{\alpha,+}(\lambda)g\in L^2(X_0),
\ g\in\Cinf_c(\sphere_{a,\reg}(\sqrt{\lambda-\ep_\alpha})),\ p<\mu-1/2,
\end{equation*}
so the preceeding expression makes $\calP_{\alpha,+}(\lambda)$ well-defined
for $\mu>1$. There is some arbitrariness in the normalization of
$\calPt_{\alpha,+}(\lambda)$. The present definition is adopted because
of its connection with the asymptotic behavior of $\calPt_{\alpha,+}(\lambda)g$
at infinity, see \cite{Vasy:Scattering}.

The backward Poisson operator is defined similarly, with
\begin{equation}\label{eq:P*-F}
\calPt_{\alpha,-}(\lambda)g=\overline{c_\alpha}\psi_{\alpha}(w^a)
\int_{\sphere_a(\sqrt{\lambda-\ep_\alpha})} e^{iw_a\cdot\omega_a}
g\,d\omega_a,\quad g\in\Cinf_c(\sphere_{a,\reg}(\sqrt{\lambda-\ep_\alpha})),
\end{equation}
\begin{equation*}\begin{split}
\calP_{\alpha,-}(\lambda)g&=\calPt_{\alpha,-}(\lambda)g
-R(\lambda-i0)((H-\lambda)
\calPt_{\alpha,-}(\lambda)g)\\
&=\calPt_{\alpha,-}(\lambda)g-R(\lambda-i0)I_a
\calPt_{\alpha,-}(\lambda)g.
\end{split}\end{equation*}

The scattering matrix relates the forward and backward Poisson operators,
i.e.\ connects incoming and outgoing data. Here we only need an expression
connecting the S-matrices to the `Green pairing'.

\begin{prop}\label{prop:boundary-pair}
Let $\alpha$ and $\beta$ be channels associated to the clusters $a$ and $b$
respectively, and suppose that $\lambda\nin\Lambda'$.
Let $u_+=\calP_{\alpha,+}(\lambda)g_+$, $u_-=\calPt_{\beta,-}(\lambda)g_-$,
$g_+\in\Cinf_c(\sphere_{a,\reg}(\sqrt{\lambda-\ep_\alpha}))$,
$g_-\in\Cinf_c(C_{b,\reg}(\sqrt{\lambda-\ep_\beta}))$. Then
\begin{equation}\begin{split}\label{eq:boundary-pair}
\langle u_+,I_b u_-\rangle&=\langle u_+,(H-\lambda)u_-\rangle
-\langle(H-\lambda) u_+,u_-\rangle\\
&=2i\sqrt{\lambda-\ep_\beta}(\langle \calS_{\alpha\beta+}(\lambda)g_+,g_-\rangle-
\delta_{\alpha\beta}\langle g_+,\calSt_{\beta\beta-}(\lambda)g_-\rangle)\\
&=2i\sqrt{\lambda-\ep_\beta}
\langle (\calS_{\alpha\beta+}(\lambda)-\delta_{\alpha\beta}
\calS_+(\lambda))g_+,g_-\rangle.
\end{split}\end{equation}
where the $L^2$ pairings on the spheres are with respect to the standard
measures normalized to have the volume of the unit sphere,
$\delta_{\alpha\beta}$ is the Kronecker delta function,
and $\calS_{\beta\beta-}(\lambda)$ is the free scattering matrix on
$X_a$ at energy $\lambda-\ep_\beta$, hence
it is a constant multiple of pull-back by the antipodal map on
$\sphere_b(\sqrt{\lambda-\ep_\beta})$.
\end{prop}

\begin{proof}
In each of the two relevant microlocal regions, namely incoming and outgoing,
one of the two functions $u_+$ and $u_-$ has trivial asymptotics. Thus,
the calculation of \cite[Section~3]{Vasy:Scattering} applies separately
in each region.
\end{proof}

With the current normalization, the S-matrix is geometric, i.e.\ under
the free evolution particles incoming at direction $\omega$ exit in
the opposite direction $-\omega$. We now introduce the relative S-matrix
(relative to free motion) as follows.
Let $p^*$ denote pull-back by the antipodal map, and let
\begin{equation}
\calS^\sharp_{\alpha\beta+}(\lambda)=\frac{1}{2i\sqrt{\lambda-\ep_\beta}}
\left(
\calS_{\alpha\beta+}(\lambda)-\delta_{\alpha\beta}
\calSt_+(\lambda)\right)p^*.
\end{equation}
If $\alpha$, $\beta$ are 2-clusters, and $V_b\in S^{-s}(X^b)$, $s>0$,
is a symbol, then the kernel of
$\calS^\sharp_{\alpha\beta+}(\lambda)$ is conormal to the diagonal
(in the sense that it is smooth for $\alpha\neq\beta$, conormal for
$\alpha=\beta$), and
if $V_b$ is Schwartz, $\calS^\sharp_{\alpha\beta+}(\lambda)$ is
a smoothing operator, i.e.\ it has a smooth kernel, as was proved
by Skibsted \cite{Skibsted:Smoothness}. If $V_b\in\Sch(X^b)$ for all $b$,
this can be seen
from \eqref{eq:boundary-pair}, for $I_a\calPt_{\alpha,-}(\lambda):
\dist(\sphere_a(\sqrt{\lambda-\ep_\alpha}))\to\Sch(X_0)$. In general,
one needs to construct a better approximation for $\calP_{\alpha+}(\lambda)$
(better than $\calPt_{\alpha+}(\lambda)$); this is what Skibsted did
in \cite{Skibsted:Smoothness}.
If $\langle w^b\rangle^\mu V_b\in L^\infty(X^b)$,
with $\mu>\dim X_a$, we
may take $g_\pm$ to be delta distributions directly (without using
Skibsted's construction, hence without a symbolic assumption),
$g_+=\delta_\omega$,
$g_-=\delta_{\omega'}$, $\omega\in \sphere_a(\sqrt{\lambda-\ep_\alpha})$,
$\omega'\in\sphere_b(\sqrt{\lambda-\ep_\beta})$.
Writing
\begin{equation}
U_\rho=(\Id-R(\lambda+i0)I_a)u^0_\rho,\ \lambda=\rho^2+\ep_\alpha,
\end{equation}
we thus deduce the following.

\begin{cor}\label{cor:pair-8}
For $\lambda\nin\Lambda'$, $\langle w^b\rangle^\mu V_b\in L^\infty(X^b)$,
$\mu>\dim X_a$,
\begin{equation}\label{eq:pair-8}
\calS^\sharp_{\alpha\beta+}
(\lambda,\omega,\omega')=\int_{\Rn} I_b U_\omega
\overline{e^{iw_b\cdot\omega'}\psi_\beta}=\int_{\Rn} I_b U_\omega
e^{-iw_b\cdot\omega'}\overline{\psi_\beta},
\end{equation}
hence $\calS^\sharp_{\alpha\beta+}$ has a continuous kernel.
\end{cor}

\section{Exponential eigenfunctions for three-body Hamiltonians}
In this section we construct exponential solutions of $(H-\lambda)u=0$
in the three-body setting. First,
for $\rho\in\Cx(X_a)$, i.e.\ $\re\rho,\im\rho\in X_a$, let
\begin{equation*}
u^0_\rho=u^0_{\alpha,\rho}=e^{i\rho\cdot w_a}\psi_\alpha(w^a).
\end{equation*}
Thus, $u^0_\rho$ is an `exponential eigenfunction' of $H_a$, namely
\begin{equation*}
(H_a-\lambda)u^0_\rho=0,\qquad \rho\cdot\rho=\lambda-\ep_\alpha.
\end{equation*}
We assume everywhere that $\dim X_a\geq 2$.

For the Hamiltonian $H$, we then seek exponential solutions $u$ of the form
\begin{equation}
u=u_\rho
=e^{i\rho\cdot w_a}(\psi_\alpha(w^a)+v),\ \rho\cdot\rho=\lambda-\ep_\alpha,
\ \rho\in \Cx(X_a),
\end{equation}
where $v$ is supposed to be `small', and $\Cx(X_a)$ denotes the
complexification of $X_a$, i.e.\ $\re\rho,\im\rho\in X_a$.
Substituting into $(H-\lambda)u=0$, we obtain
\begin{equation}
(\Delta+2\rho\cdot D_w+V_a+I_a-\ep_\alpha)v=-I_a\psi_\alpha.
\end{equation}
The right hand side decays at infinity since $\psi_\alpha$
does so in $X^a$, and $I_a$ decays away from $\cup_{\#b=2,\ b\neq a}C_b$.
More precisely, we have the following lemma.

\begin{lemma}
Suppose that $\mu>0$. There exists $C>0$ with the following property.
If $\langle w^b\rangle^{\mu}V_b\in L^\infty(X^b)$ for all
two-clusters $b$ with $b\neq a$ then
\begin{equation}\label{eq:I_a-w_a}
|I_a|\leq C\langle w_a\rangle^{-\mu}
\langle w^a\rangle^{\mu}\sup|\langle w^b\rangle^{\mu}V_b|
\end{equation}
\end{lemma}

\begin{proof}
As $X_a\cap X_b=\{0\}$ for $b\neq a$, $X^a\oplus
X^b=X_0$, hence for some $C'>0$
\begin{equation*}
\langle w^a\rangle\langle w^b\rangle\geq C'\langle w\rangle.
\end{equation*}
Thus, for $\mu>0$,
\begin{equation}
\langle w^b\rangle^{-\mu}\leq C \langle w_a\rangle^{-\mu}
\langle w^a\rangle^{\mu},
\end{equation}
proving the lemma.
\end{proof}

\begin{cor}
Let $\mu>0$ and $V_b$ as above. Then
\begin{equation}\label{eq:I_a-psi_a}
I_a\psi_\alpha\in L^2_p(X_0),\ p<\mu-\frac{\dim X_a}{2}.
\end{equation}
\end{cor}

Thus, we need to construct a right inverse $G(\rho)$ to
\begin{equation*}
P(\rho)=\Delta+2\rho\cdot D_w+V-\ep_\alpha
\end{equation*}
that can be applied to elements of $L^2_p(X_0)$.
Once this is done,
\begin{equation*}
u(\rho)=e^{i\rho\cdot w_a}(\psi_\alpha(w^a)-G(\rho) I_a\psi_\alpha)
\end{equation*}
is the solution to the original problem.
Below we write
\begin{equation}\begin{split}\label{eq:P_a-def}
&P_0(\rho)=\Delta+2\rho\cdot D_w-\ep_\alpha,\\
&P_a(\rho)=\Delta+2\rho\cdot D_w+V_a-\ep_\alpha\\
&\qquad=(\Delta_{w_a}+2\rho\cdot D_{w_a}-\ep_\alpha)
+(\Delta_{w^a}+2\rho\cdot D_{w^a}+V_a(w^a)).
\end{split}\end{equation}
Since a right inverse $G_a(\rho)$ of $P_a(\rho)$ can be constructed
explicitly, perturbation theory will give the existence of $G(\rho)$,
provided that $I_a$ is small.

It is convenient to represent $\rho$ as
\begin{equation*}
\rho=z\nu+\rho_\perp,\ \rho_\perp,\nu\in X_a,\ z\in\Cx,\ \rho_\perp\cdot\nu=0.
\end{equation*}
We will take $|\rho_\perp|$ sufficiently small.
To see why the
size of $\rho_\perp$ matters, consider
\begin{equation*}
G_0(\rho)=\Frinv (|\xi|^2+2\rho\cdot\xi-\ep_\alpha)^{-1} \Fr,
\end{equation*}
so $P_0(\rho)G_0(\rho)=\Id$ e.g.\ on Schwartz functions.
Thus, on the Fourier transform side
$G_0(\rho)$ acts via multiplication by
$(|\xi|^2+2\rho\cdot\xi-\ep_\alpha)^{-1}$. For $z\nin\Real$, this distribution
is conormal to
\begin{equation}\begin{split}
S(\rho)&=\{\xi\in X_0:\ |\xi|^2+2\re\rho\cdot\xi-\ep_\alpha=0,
\ \im\rho\cdot\xi=0\}\\
&=\{\xi\in X_0:\ (\xi+\rho_\perp)^2=\rho_\perp^2+\ep_\alpha,
\ \nu\cdot\xi=0\}.
\end{split}\end{equation}
Note that $S(\rho)$ actually depends only on $\rho_\perp$ and $\nu$,
not on $z$. Now, for $|\rho_\perp|<\sqrt{-\ep_\alpha}$
(note that $\ep_\alpha<0$), $S(\rho)=\emptyset$, so $P_0(\rho)$ is elliptic
`at infinity' in a sense discussed by Melrose \cite{RBMSpec}, namely
as an element of $\Diffsc^2(\overline{X_0})$, $\overline{X_0}$ being
the radial compactification of $X_0$.

Below we assume that
\begin{equation}\label{eq:rho_perp-bd}
|\rho_\perp|<\sqrt{-\ep_\alpha},
\end{equation}
where $P_0(\rho)$ is elliptic.
This does {\em not} mean that $P(\rho)$
itself is elliptic; indeed $P_a(\rho)$ cannot be such thanks to the
bound state $\psi_\alpha$. For
\begin{equation}\label{eq:ep_1-def}
\evpr\in [|\rho_\perp|^2+\ep_\alpha,0)\setminus
\Lambda'_a,
\end{equation}
let $e_a(\evpr)$ be the orthogonal projection
to the $L^2$ eigenfunctions of $H^a$
with eigenvalue $\leq\evpr$,
\begin{equation*}
e_a(\evpr)=\sum_{\alpha':\ \ep_{\alpha'}\leq \evpr}
(\psi_{\alpha'}\otimes\overline{\psi_{\alpha'}})\in\bop(L^2(X^a),L^2(X^a)),
\end{equation*}
and let $E_a$ be its extension to $X_0$ via
tensoring by $\Id_{X_a}$, so
\begin{equation}
E_a(\evpr)=\sum_{\alpha'}
\Id_{X_a}\otimes(\psi_{\alpha'}\otimes\overline{\psi_{\alpha'}}).
\end{equation}
Since eigenvalues of $H^a$ can only accumulate at $\Lambda_a=\{0\}$,
$e_a$ is finite rank. We also let
\begin{equation}\label{eq:ep_0-def}
\evth=\inf\left(\Lambda'_a\cap(\evpr,+\infty)\right)>\evpr.
\end{equation}
{\em The particular choice of $\evpr$, provided that it is sufficiently
close to $0$, does not play a major role in our arguments, so we usually
simply write $e_a$ for $e_a(\evpr)$, etc.}

We restrict the region \eqref{eq:rho_perp-bd}
slightly further and work in the region
\begin{equation}
\Cx(X_a)_\alpha^\circ
=\{(z,\nu,\rho_\perp):\ \im z\neq 0,
\ |\rho_\perp|^2+\ep_\alpha\in(\ep_\alpha,\evth)\setminus\Lambda'_a\},
\end{equation}
i.e.\ we also assume that $|\rho_\perp|^2+\ep_\alpha$ is not an eigenvalue
of $H^a$. Again, we do not indicate $\evth$ explicitly in the notation.

Since the ranges of $E_a$ and $\Id-E_a$ play a rather different role
below, we introduce weighted spaces that reflect this. So for $p\in\Real$
we let
\begin{equation}\begin{split}\label{eq:calH_p-def}
\calH_p=&(L^2_p(X_a)\otimes \Range e_a)\oplus (L^2(X_a)\otimes
\Range (\Id-e_a))\\
&\qquad\subset (L^2_p(X_a)\otimes \Range e_a)\oplus L^2(X_0),
\end{split}\end{equation}
with $e_a$ considered as a bounded operator on $L^2(X^a)$.
Thus, we allow weights on the range of $E_a$, but not on its orthocomplement.
Again, $\Range e_a$ is finite dimensional, hence it is closed
in $L^2_r(X^a)$ for all $r\in\Real$, while $\Range(\Id-e_a)$ is
closed in $L^2(X^a)$, so $L^2(X_a)\otimes \Range (\Id-e_a)$ is a closed
subspace of $L^2(X_0)$. Thus, for all $p\in\Real$, $\calH_p$ is a Hilbert
space with norms induced on the summands by the $L^2_p(X_a)$ and
$L^2(X_0)$ norms respectively.

We start the construction of $G(\rho)$ by analyzing $G_a(\rho)$.
In view of \eqref{eq:P_a-def},
taking the Fourier transform in $X_a$ makes the invertibility of
$P_a(\rho)$ into a question on the behavior of the resolvent
of $H^a=\Delta_{w^a}+V_a(w^a)$, uniformly across the spectrum.
That is,
\begin{equation*}
\Fr_{X_a} P_a(\rho)\Fr_{X_a}^{-1}=
\Delta_{X^a}+V_a-(\ep_\alpha-|\xi_a|^2-2\rho\cdot\xi_a),
\end{equation*}
acting pointwise in $\xi_a$, so
\begin{equation*}
\Fr_{X_a} G_a(\rho)\Fr_{X_a}^{-1}=R^a(\ep_\alpha-|\xi_a|^2-2\rho\cdot\xi_a),
\end{equation*}
with $R^a(\sigma)=(H^a-\sigma)^{-1}$, provided we show that this makes sense
-- the only issues being the behavior for real $\sigma$ and bounds as
$|\xi_a|\to\infty$.

So let
\begin{equation*}
F_\rho:\xi_a\mapsto \ep_\alpha-|\xi_a|^2-2\rho\cdot\xi_a.
\end{equation*}
Thus,
\begin{equation}\begin{split}\label{eq:re-z-im-z}
&\im F_\rho(\xi_a)=-2(\im z)(\nu\cdot\xi_a)\\
&\re F_\rho(\xi_a)=\ep_\alpha+\rho_\perp^2
-((\xi_a)_\perp+\rho_\perp)^2
-(\nu\cdot\xi_a)^2-2(\re z)(\nu\cdot\xi_a)\\
&\qquad\qquad\qquad
\leq \ep_\alpha+\rho_\perp^2-(\nu\cdot\xi_a)^2-2(\re z)(\nu\cdot\xi_a)\\
&\qquad\qquad\qquad\leq \ep_\alpha+\rho_\perp^2+(\re z)^2,
\end{split}\end{equation}
where $(\xi_a)_\perp$ is the orthogonal projection of $\xi_a$ to the
orthocomplement of the span of $\nu$, and
$\nu\cdot\xi_a$ is the component of $\xi_a$ parallel to $\nu$.
If $F_\rho(\xi_a)$ is real and $\im z\neq 0$, then
$\nu\cdot\xi_a=0$, hence
\begin{equation*}
\ep_\alpha-|\xi_a|^2-2\rho\cdot\xi_a=\ep_\alpha+\rho_\perp^2
-(\xi_a+\rho_\perp)^2\leq\ep_\alpha+\rho_\perp^2<0
\end{equation*}
under our assumptions. In fact,
\begin{equation*}
\re F_\rho(\xi_a)\leq \ep_\alpha+\rho_\perp^2
-2(\re z)(\nu\cdot\xi_a)
=\ep_\alpha+\rho_\perp^2+2\frac{\re z}{\im z}\im F_\rho(\xi_a),
\end{equation*}
so $\re F_\rho(\xi_a)>\ep_\alpha+\rho_\perp^2$, which holds in particular
if $\re F_\rho(\xi_a)\geq 0$, implies that
$\nu\cdot\xi_a$ is non-zero and has the same
sign as $-\re z$. Indeed, we deduce that
\begin{equation}\label{eq:sgn-im-F}
\re F_\rho(\xi_a)>\ep_\alpha+\rho_\perp^2\Rightarrow
\frac{\re z}{\im z}\im F_\rho(\xi_a)>0.
\end{equation}
In fact, we deduce the quantitative bound
\begin{equation}\label{eq:sgn-im-F-p}
\re F_\rho(\xi_a)>\ep_\alpha+\rho_\perp^2\Rightarrow
\im F_\rho(\xi_a)>(\re F_\rho(\xi_a)-(\ep_\alpha+\rho_\perp^2))
\frac{|\im z|}{|\re z|}
\end{equation}

Moreover, if $\dim X_a\geq 2$, as is assumed throughout this paper,
$F_\rho:X_a^*\to\Cx$ (the latter considered as a 2-dimensional
real manifold)
has a surjective differential unless $(\xi_a)_\perp=-\rho_\perp$.
Indeed, since $d\im F_\rho$ is nonzero, and is a multiple of
$d(\nu\cdot\xi_a)$,
$d\im F_\rho$ and $d\re F_\rho$ are linearly independent if and only if
$d(((\xi_a)_\perp+\rho_\perp)^2)\neq 0$, i.e.\ if and only if
$(\xi_a)_\perp\neq-\rho_\perp$.

Note also that for any fixed $\rho$ there exists $C>0$ such that
\begin{equation}\label{eq:re F<0}
\re F_\rho(\xi_a)<C-|\xi_a|^2/2.
\end{equation}

The structure of $R^a(\sigma)$
corresponds to that of $R^a_0(\sigma)$ and
$(\Id+V_a R^a_0(\sigma))^{-1}$, where
$R^a_0(\sigma)$ denotes the resolvent of $\Delta_{X^a}$, for
$R^a(\sigma)=R^a_0(\sigma)(\Id+V_a R^a_0(\sigma))^{-1}$. It satisfies
automatically that for $\re\sigma<\inf\spec H^a$,
\begin{equation*}
\|R^a(\sigma)\|_{\bop(L^2,L^2)}\leq (\inf\spec H^a-\re\sigma)^{-1}.
\end{equation*}
In view of \eqref{eq:re F<0}, for $\phi\in\Cinf_c(X_a^*)$ identically $1$
on a large enough ball, $G_a\Frinv_{X_a}(1-\phi)\Fr_{X_a}$ is bounded
on $L^2$. Hence, we only need to be concerned about what happens in
a compact set in $X_a^*$. We also mention two other bounds that hold
by the selfadjointness of $H^a$ and its spectral properties, namely
\begin{equation}\begin{split}\label{eq:spec-bounds}
&\|(\Id-e_a)R^a(\sigma)\|_{\bop(L^2,L^2)}\leq |\re\sigma-\evth|^{-1},
\ \re\sigma<\evth,\\
&\|R^a(\sigma)\|_{\bop(L^2,L^2)}\leq |\im\sigma|^{-1},\ \im\sigma\neq 0,
\end{split}\end{equation}
$\evth$ as in \eqref{eq:ep_0-def}.

Now $R^a(\sigma)$ is analytic in $\sigma$ for
$\im\sigma\neq 0$, with values in
bounded operators
$L^2(X^a)\to H^2(X^a)$,
and from $\im\sigma>0$ (and from
$\im\sigma<0$) it extends to be smooth to $\Cx\setminus[0,+\infty)$
away from the eigenvalues
of $H^a$, where it has a simple pole.
With $e_a$ denoting the projection to the $L^2$ eigenspace of $H^a$
with eigenvalues $\leq \evpr$, $\evpr\in[\rho_\perp^2+\ep_\alpha,0)$,
and $\evth$ given by \eqref{eq:ep_0-def}
as before, $R^a(\sigma)(\Id-e_a)$ is smooth on
$\Cx\setminus(\evth,+\infty)$.
On the range of $e_{a,\ep}$, the projection
to the eigenspace with eigenvalue $\ep$, $R^a(\sigma)$ is multiplication
by $(\ep-\sigma)^{-1}$. This is a locally integrable function of $\sigma$
near $\ep$ (in $\Cx$!), so the application of $R^a(F_\rho(.))$ to
$\hat u=\Fr_{X_a}u$ is well defined, provided that at every point $\xi_a$
with $F_\rho(\xi_a)=\ep$, the differential of $F_\rho$ is surjective,
i.e.\ $F_\rho(\xi_a)=\ep$ implies that $(\xi_a)_\perp\neq-\rho_\perp$.
But if $(\xi_a)_\perp=-\rho_\perp$ and $F_\rho(\xi_a)$ is real,
then $\nu\cdot\xi_a=0$, hence $F_\rho(\xi_a)=\ep_\alpha+\rho_\perp^2$.
Now let
\begin{equation*}\begin{split}
G_a(\rho)&=\Frinv_{X_a}R^a(\ep_\alpha-|.|^2-2\rho\cdot .)\Fr_{X_a}\\
&=\Frinv_{X_a}R^a(\ep_\alpha-|.|^2-2\rho\cdot .)(\Id-E_a)\Fr_{X_a}
+\Frinv_{X_a}R^a(\ep_\alpha-|.|^2-2\rho\cdot .)E_a\Fr_{X_a};
\end{split}\end{equation*}
both terms are well defined by the preceeding considerations when applied
to functions in $\Sch(X_a;L^2(X^a))$. Indeed, application of $G_a(\rho)$
to the range of $E_{a,\ep}$ is the only issue, and there,
with $u=v\otimes \psi_{\alpha'}$,
$v\in L^2_p(X_a)$, $\hat v=\Fr_{X_a}v$,
\begin{equation*}
(\Fr_{X_a}E_{a,\ep}G_a(\rho)v)(\xi_a,.)=(\ep-F_\rho(\xi_a))^{-1}\hat v
\otimes \psi_{\alpha'},
\end{equation*}
so the mapping properties of $G_a(\rho)$ on $\Range E_a$ are given by
the two-body results of Weder \cite{Weder:Generalized}. In particular
$E_a G_a(\rho)$ is
well defined for functions in $L^2_p(X_a)\otimes \Range e_a$, $p>0$.
Hence we deduce the following result.

\begin{prop}
Suppose that $(z,\nu,\rho_\perp)\in\Cx(X_a)_\alpha^\circ$.
The operator
\begin{equation*}
G_a(\rho)=\Frinv_{X_a}R^a(\ep_\alpha-|.|^2-2\rho\cdot .)\Fr_{X_a}
\end{equation*}
is a bounded operator
$\calH_p\to \calH_r$ for
$p>0$, $r<0$, $r<p-1$. It satisfies
\begin{equation}\begin{split}\label{eq:P_aG_a-Id}
&P_a(\rho)G_a(\rho)=\Id:\calH_p\to \calH_p,\\
&G_a(\rho)P_a(\rho)=\Id:\calH_p\to \calH_p.
\end{split}\end{equation}
It is continuous in $\rho\in\Cx(X_a)_\alpha^\circ$ and
analytic in $z\in\Cx\setminus\Real$.
Moreover, for $p>0$, $r<0$, $r<p-1$, $\rho_\perp$ fixed,
for any $C>0$, $G_a(\rho)$ is uniformly bounded in $\bop(\calH_p,\calH_r)$
in $|\im z|\geq C|\re z|$,
and $\slim_{|z|\to\infty} G_a(\rho)=0$
as an operator in $\bop(\calH_p,\calH_r)$,
provided that $|z|\to\infty$ in the region
$|\im z|\geq C|\re z|$.
\end{prop}

\begin{proof}
All of the claims follow from the previous argument, except the
behavior of $G_a(\rho)$ as $\rho\to\infty$. That in turn follows
from
\begin{equation}\begin{split}
&\lim_{\rho\to\infty}\|E_a G_a(\rho)\|_{\bop(\calH_p,\calH_r)}
= 0,\ p>0,r<0,
r<p-1,\\
&\|(\Id-E_a)G_a(\rho)\|_{\bop(\calH_p,\calH_p)}\leq C,\\
&(\Id-E_a)G_a(\rho)\to 0\ \text{strongly on}\ \calH_p.
\end{split}\end{equation}
The first estimate here is a two-body result, relying on
a similar estimate for $\|E_{a,\ep} G_a(\rho)\|$ for each $\ep\in\pspec(H^a)$,
see \cite{Weder:Generalized}. Indeed, with $u=v\otimes \psi_{\alpha'}$,
$v\in L^2_p(X_a)$, $\hat v=\Fr_{X_a}v$,
\begin{equation*}
(\Fr_{X_a}E_{a,\ep}G_a(\rho)v)(\xi_a,.)=(\ep-F_\rho(\xi_a))^{-1}\hat v
\otimes \psi_{\alpha'}.
\end{equation*}
The uniform estimate for $(\Id-E_a)G_a(\rho)$ holds because
$R^a(\sigma)(\Id-e_a)$ is uniformly bounded as an operator on $L^2(X^a)$
as long as $\sigma$ is uniformly bounded away from $[0,+\infty)$, which
holds for $\sigma=F_\rho(\xi_a)$
provided $|\im z|\geq C|\re z|$ by \eqref{eq:sgn-im-F-p}. Then by the
Parseval's formula, and with $\hat u=\Fr_{X_a}u$,
\begin{equation*}\begin{split}
\|G_a(\rho)(\Id-E_a)u&\|_{L^2(X_0)}^2
=(2\pi)^{-n}
\int_{X_a}\|R^a(F_\rho(\xi_a))\hat u(\xi_a,.)\|^2_{L^2(X^a)}\,d\xi_a\\
&\leq (2\pi)^{-n}M^2\int_{X_a}\|\hat u(\xi_a,.)\|^2_{L^2(X^a)}\,d\xi_a
=(2\pi)^{-n}M^2\|u\|^2_{L^2(X_0)},\\
&\qquad M=\sup\{\|R^a(F_\rho(\xi_a))(\Id-e_a)\|_{\bop(L^2(X^a),L^2(X^a))} :
\ \xi_a\in X_a\}.
\end{split}\end{equation*}
Moreover, as $\im z\to\infty$, $\im F_\rho(\xi_a)\to \infty$ for almost every
$\xi_a$, namely for $\xi_a$ such that $\xi_a\cdot\nu\neq0$. But
$R^a(\sigma)\to 0$ as a bounded operator on $L^2(X^a)$ as $\im\sigma\to\infty$
by \eqref{eq:spec-bounds}. Since
\begin{equation*}
\|R^a(F_\rho(\xi_a))\hat u(\xi_a,.)\|^2
\leq M^2\|\hat u(\xi_a,.)\|^2_{L^2(X^a)},
\end{equation*}
and $\|\hat u(\xi_a,.)\|^2_{L^2(X^a)}\in L^1(X_a)$,
the dominated convergence theorem implies that
$G_a(\rho)(\Id-E_a)\to 0$ strongly.
\end{proof}

Since
\begin{equation*}
P(\rho)G_a(\rho)v=(\Id+I_a G_a(\rho))v,\ v\in\calH_p,\ p>0,
\end{equation*}
we next investigate $I_a$ on $\calH_r$.

\begin{lemma}\label{lemma:I_a-norm}
Suppose that $\mu>0$, and $\langle w^b\rangle^\mu
V_b\in L^\infty(X^b)$ for all $b\neq a$.
The multiplication operator $I_a$ is in $\bop(\calH_r,\calH_p)$ provided
that $p\leq r+\mu$, $p\leq \mu$, $r\geq -\mu$.

Moreover, there exists $C>0$
(independent of $V_b$) such that the norm of $I_a$
as such an operator is bounded by $C\max_b\sup(\langle w^b\rangle^\mu |V_b|)$.
\end{lemma}

\begin{proof}
We decompose $I_a$ as a matrix corresponding to the direct sum in
\eqref{eq:calH_p-def}. Since $I_a$ is bounded on $L^2(X_0)$,
it follows that $(\Id-E_a)I_a(\Id-E_a):\calH_r\to\calH_p$ for all $r$ and
$p$. Moreover, by \eqref{eq:I_a-w_a}, for all $\alpha'$, $\alpha''$,
\begin{equation*}
\langle w_a\rangle^\mu
\int_{X^a} I_a \overline{\psi_{\alpha'}}\psi_{\alpha''}\,dw_a
\end{equation*}
is bounded on $X_a$, hence $E_a I_a E_a:\calH_r\to\calH_p$ for all $r$ and
$p$ with $p\leq r+\mu$. In addition,
$\langle w_a \rangle^\mu I_a\overline{\psi_{\alpha'}}\in L^\infty(X_0)$,
so $E_a I_a:L^2(X_0)\to \calH_p$ for $p\leq\mu$. Similarly,
$\langle w_a\rangle^\mu\psi_{\alpha'}I_a\in L^\infty(X_a;L^2(X^a))$,
so $I_a E_a:\calH_r\to L^2(X_0)$ provided $r\geq -\mu$. As $\Id-E_a$ is
a bounded operator on $L^2(X_0)$, this shows that $I_a:\calH_r\to\calH_p$
as stated.
\end{proof}

Combining the preceeding proposition and lemma we deduce the following.

\begin{cor}
Suppose that $\mu>1$, and $\langle w^b\rangle^\mu
V_b\in L^\infty(X^b)$ for all $b\neq a$. Let $p$ satisfy $0<p<\mu$.
Then $I_a G_a(\rho)\in\bop(\calH_p,\calH_p)$, continuous in $\rho\in
\Cx(X_a)_\alpha^\circ$, and analytic in $z$ in this region.
Moreover, for $\rho_\perp$ fixed and
$C>0$, $I_a G_a(\rho)$ is uniformly bounded in
$|\im z|\geq C|\re z|$, $|z|>1$, and $\slim_{|z|\to\infty}I_aG_a(\rho)=0$
inside this region.
\end{cor}

We also need to consider the invertibility properties of $\Id+I_a G_a(\rho)$
on $\calH_p$. As $I_a G_a(\rho)\in\bop(\calH_p,\calH_p)$, $\Id+I_a G_a(\rho)$
is invertible and $\|(\Id+I_a G_a(\rho))^{-1}\|_{\bop(\calH_p,\calH_p)}<2$
provided that $\|I_a G_a(\rho)\|_{\bop(\calH_p,\calH_p)}<1/2$. In view
of the uniform boundedness of $G_a(\rho)$ as a map in $\bop(\calH_p,\calH_p)$,
there exists $\delta'>0$ such that if
$\|I_a\|_{\bop(\calH_r,\calH_p)}<\delta'$ then
$\|I_a G_a(\rho)\|_{\bop(\calH_p,\calH_p)}<1/2$. Combining this with
the norm estimate of Lemma~\ref{lemma:I_a-norm} leads to the following
theorem.

\begin{thm}\label{thm:G(rho)-exist)}
Suppose that $\rho_\perp$ satisfies $|\rho_\perp|^2+\ep_\alpha\in
(\ep_\alpha,\evth)\setminus\Lambda'_a$,
$\mu>\max(p,1)$, $p>0$, $r<0$, $r<p-1$.
There exists $\delta>0$ with the following property. Suppose that
for all $b\neq a$, $\sup |\langle
w^b\rangle^\mu V_b|<\delta$.
Then the operator
\begin{equation}\label{eq:G(rho)-def}
G(\rho)=G_a(\rho)(\Id+I_a G_a(\rho))^{-1}:\calH_p\to \calH_r
\end{equation}
satisfies
\begin{equation}\begin{split}\label{eq:PG-Id}
&P(\rho)G(\rho)=\Id:\calH_p\to \calH_p,\\
&G(\rho)P(\rho)=\Id:\calH_p\to \calH_p.
\end{split}\end{equation}
Moreover, $G(\rho)$ is a continuous function of $\rho$ in
$\Cx(X_a)_\alpha^\circ$, and an analytic function of $z\in\Cx\setminus\Real$,
and $\slim_{|z|\to\infty}G(\rho)=0$ as a map $\calH_p\to\calH_r$ provided that
$|z|\to\infty$ in $|\im z|\geq C|\re z|$, $C>0$.
\end{thm}

\begin{proof}
With $G(\rho)$ as in \eqref{eq:G(rho)-def}, $G(\rho):\calH_p\to \calH_r$,
$p>0$, $r<0$, $r<p-1$, since $G_a(\rho)$ has these mapping properties,
and $(\Id+I_a G_a(\rho))^{-1}$ is bounded on $\calH_p$.
Now,
\begin{equation}\begin{split}
&P(\rho)G(\rho)=(P_a(\rho)+I_a)G_a(\rho)(\Id+I_a G_a(\rho))^{-1},\\
&(P_a(\rho)+I_a)G_a(\rho)=\Id+I_a G_a(\rho):\calH_p\to \calH_p,
\end{split}\end{equation}
proving the first line of \eqref{eq:PG-Id}. The second line follows
from the identity
\begin{equation*}
(\Id+I_a G_a(\rho))^{-1}=\Id-(\Id+I_a G_a(\rho))^{-1}I_a G_a(\rho),
\end{equation*}
and $G_a(\rho)P_a(\rho)=\Id$ on $\calH_p$, $p>0$. The limiting behavior
follows from the uniform boundedness of $(\Id+I_a G_a(\rho))^{-1}$
on $\calH_p$, and from
\begin{equation*}
G(\rho)=G_a(\rho)-G_a(\rho)(\Id+I_a G_a(\rho))^{-1}(I_a G_a(\rho))
\end{equation*}
with the last factor tending to $0$ strongly on $\calH_p$, and the
other factors remaining bounded.
\end{proof}

Since $I_a\psi_\alpha\in \calH_p$ for some $p>0$ if $\mu>\dim X_a/2$,
due the \eqref{eq:I_a-psi_a}, we deduce the following corollary.

\begin{cor}\label{cor:u_rho-exists}
Suppose that $\mu>\dim X_a/2$,
$|\rho_\perp|^2+\ep_\alpha\in(\ep_\alpha,\evth)\setminus\Lambda'_a$.
Then
\begin{equation}
u_\rho=u^0_\rho-e^{i\rho\cdot w_a}G(\rho)I_a \psi_\alpha=e^{i\rho\cdot w_a}
(\psi_\alpha-G(\rho)I_a \psi_\alpha)
\end{equation}
satisfies $P(\rho)u_\rho=0$, and
\begin{equation}
u^0_\rho=e^{i\rho\cdot w_a}(\Id+G_a(\rho)I_a)e^{-i\rho\cdot w_a}u_\rho.
\end{equation}
Moreover, $u_\rho-u^0_\rho\to 0$ in $\calH_r$, $r<0$,
$r<\mu-1-\frac{\dim X_a}{2}$,
as $|z|\to\infty$ in $|\im z|>C|\re z|$, $C>0$.
\end{cor}

\section{Limit as $\rho$ goes to the reals}
As $\rho$ becomes real, the structure of the operator
$G_a(\rho)$ degenerates since $|\xi|^2+2\rho\cdot\xi-\ep_\alpha$
becomes real. Many of the details of the following calculations
are similar to the corresponding two-body calculations, see
e.g.\ Weder \cite{Weder:Generalized}.
We proceed as follows.

Recall that for $\rho\in\Cx(X_a)$,
$\im\rho\neq 0$, there exist unique $\nu\in X_a$, $z\in \Cx$,
$\rho_\perp\in X_a$ such that
\begin{equation*}
\rho=z\nu+\rho_\perp,\quad|\nu|^2=1,\ \rho_\perp\cdot\nu=0,\ \im z>0.
\end{equation*}
Alternatively, the inequality $\im z>0$ can be replaced by $\im z<0$.
Here we keep $\im z>0$ for the sake of definiteness.
The behavior as $\rho$ approaches $X_a$ then corresponds to $z$ approaching
the real axis. Now, $u_\rho$ solves $(H-\lambda)u_\rho=0$,
$\lambda=\rho^2+\ep_\alpha$, i.e.\ has energy $\rho^2+\ep_\alpha$.
Thus, the form of the limit as $\rho\to X_a$ depends on the nature of
the spectrum of $H$ near $\lambda=\rho^2+\ep_\alpha$. Here we restrict
ourselves to $|\rho_\perp|^2<-\ep_\alpha$. Note that with $\evth$
as in \eqref{eq:ep_0-def},
\begin{equation}
\lambda<0\Miff\rho^2=z^2+\rho_\perp^2<-\ep_\alpha;\ \lambda<\evth\Miff
\rho^2=z^2+\rho_\perp^2<\evth-\ep_\alpha.
\end{equation}
Below we consider
\begin{equation}
\re z\geq 0,\ \im z\geq 0.
\end{equation}
and show that
$u_\rho=u_{\alpha,\rho}$ extends to a continuous
function of $(z,\nu,\rho_\perp)$ in
\begin{equation*}\begin{split}
\Cx(X_a)_\alpha^+=&\{(z,\nu,\rho_\perp):\ \im z>0,
\ |\rho_\perp|^2+\ep_\alpha\in(\ep_\alpha,\evth)
\setminus\Lambda'_a\}\\
&\cup \{(z,\nu,\rho_\perp):\ |\rho_\perp|^2+\ep_\alpha\in(\ep_\alpha,\evth)
\setminus\Lambda'_a,\ z\geq 0,\\
&\qquad\qquad\qquad\qquad
z^2+\rho_\perp^2+\ep_\alpha\in(\ep_\alpha,\evth)\setminus\Lambda'_a\}.
\end{split}\end{equation*}

\begin{rem}
Below we often consider $\rho\in\Cx(X_a)_\alpha^+$, and keep writing
$u_\rho$ even when the projection of $\rho$ to $\Cx(X_a)$ lies in $X_a$.
\end{rem}

First consider $G_a(\rho)E_a$. It is better to consider this as
$\sum_{\ep\in\pspec(H^a)}G_a(\rho)E_{a,\ep}$, where $E_{a,\ep}$ is projection
to the $\ep$ eigenspace of $H^a$ tensored with the identity map on $X_a$.
Now, for $z\nin\Real$,
\begin{equation}\begin{split}
\Fr_{X_a} G_a(\rho)&E_{a,\ep}\Frinv_{X_a} u
=(|\xi_a|^2+2\rho\cdot\xi_a-\ep_\alpha+\ep)^{-1}E_{a,\ep}u\\
&=(|\xi_a+\re\rho|^2+2i(\im z)\nu\cdot\xi_a-\re\rho^2-\ep_\alpha+\ep)^{-1}
H(\nu\cdot\xi_a)E_{a,\ep}u\\
&\quad+(|\xi_a+\re\rho|^2+2i(\im z)\nu\cdot\xi_a-\re\rho^2-\ep_\alpha+\ep)^{-1}
H(-\nu\cdot\xi_a)E_{a,\ep}u,
\end{split}\end{equation}
where $H$ is the Heaviside step function, so $H=1$ on $(0,+\infty)$,
$H=0$ on $(-\infty,0)$.
Thus, if $\im z>0$, letting $\im z\to 0$, yields
\begin{equation}\begin{split}\label{eq:G_a-lim}
\lim_{\im z\to 0}\Fr_{X_a}& G_a(\rho)E_{a,\ep}\Frinv_{X_a} u\\
&=(|\xi_a+\re\rho|^2-(\re\rho^2+\ep_\alpha-\ep-i0))^{-1}
H(\nu\cdot\xi_a)E_{a,\ep}u\\
&\qquad+(|\xi_a+\re\rho|^2-(\re\rho^2+\ep_\alpha-\ep+i0))^{-1}
H(-\nu\cdot\xi_a)E_{a,\ep}u.
\end{split}\end{equation}
This calculation makes sense for $\Frinv_{X_a} u\in\calH_s$, $s>1/2$.
Note that if $\ep>\lambda=\re\rho^2-\ep_\alpha$, then 
$|\xi_a+\re\rho|^2-(\re\rho^2+\ep_\alpha-\ep)$ does not vanish for any
$\xi_a$ (as it is positive), so in $E_{a,\ep}$ the factors with $+i0$ and
$-i0$ are the same.
Conjugating $G_a(\rho)$ by $e^{i\re\rho\cdot w_a}$ replaces
$\xi_a$ by $\xi_a+\re\rho$. Writing $H(-\nu\cdot\xi_a)=1-H(\nu\cdot\xi_a)$
and $\lambda=\re\rho^2+\ep_\alpha$,
we thus deduce that
\begin{equation}\begin{split}\label{eq:G_a-lim-16}
&\lim_{\im z\to 0}e^{i\re\rho\cdot w_a} G_a(\rho)E_{a,\ep}
e^{-i\re\rho\cdot w_a}\\
&=\Frinv
(|\xi|^2-(\lambda-\ep+i0))^{-1}\Fr E_{a,\ep}\\
&\quad +\Frinv
H(\nu\cdot(\xi_a-\re\rho))((|\xi|^2-(\lambda-\ep-i0))^{-1}
-(|\xi|^2-(\lambda-\ep+i0))^{-1})\Fr E_{a,\ep}\\
&=R_a(\lambda+i0)E_{a,\ep}+\Frinv H(\nu\cdot(\xi_a-\re\rho)
\Fr [R_a(\lambda-i0)-R_a(\lambda+i0)]E_{a,\ep}.
\end{split}\end{equation}

Now consider $G_a(\rho)(\Id-E_a)$. This is analytic in
\begin{equation}
\{z\in\Cx:\ z^2\nin[\evpr-\ep_\alpha-\rho_\perp^2,+\infty)\},
\end{equation}
hence in $\Cx(X_a)^+_\alpha$. Indeed, since
\begin{equation}
G_a(\rho)(\Id-E_a)=\Fr^{-1}_{X_a}R^a(\ep_\alpha-|\xi_a|^2-2\rho\cdot\xi_a)
(\Id-E_a)\Fr_{X_a},
\end{equation}
we need to analyze the range of the map
$F_\rho:\xi_a\mapsto \ep_\alpha-|\xi_a|^2-2\rho\cdot\xi_a$, which
has been done in \eqref{eq:re-z-im-z}-\eqref{eq:sgn-im-F}.
Thus, when $\im z\to 0$
and $0\leq \re z<\sqrt{\evpr-\ep_\alpha-\rho_\perp^2}$ implies
$\re F_\rho< \evpr$ by \eqref{eq:re-z-im-z}, and $R^a(\sigma)(\Id-e_a)$
is analytic in $\re\sigma<\evpr$.
We only need the additional observation that
as $\im z\to 0$, $\im F_\rho\to 0$.
Thus, in $\im z>0$, by the dominated convergence theorem,
\begin{equation}
\lim_{\im z\to 0}G_a(\rho)(\Id-E_a)
=\Fr^{-1}_{X_a}R^a(\ep_\alpha-(|\xi_a|^2+2\rho\cdot\xi_a-i0))
(\Id-E_a)\Fr_{X_a}.
\end{equation}
Again,
conjugating by $e^{i\rho\cdot w_a}$ replaces $\xi_a$ by $\xi_a+\rho$. Thus,
\begin{equation}\begin{split}
&\lim_{\im z\to 0}e^{i\re\rho\cdot w_a}G_a(\rho)(\Id-E_a)e^{-i\re\rho\cdot w_a}\\
&\qquad=\Fr^{-1}_{X_a}R^a(\ep_\alpha+\rho^2-(|\xi_a|^2+i0))
(\Id-E_a)\Fr_{X_a}=R_a(\lambda+i0)(\Id-E_a),
\end{split}\end{equation}
where we took into account that $\ep_\alpha+\rho^2=\lambda$. Note that
$R_a(\lambda+i0)(\Id-E_a)=R_a(\lambda-i0)(\Id-E_a)$ as $\lambda<\evth$.

Combining these two results we deduce the following proposition.

\begin{prop}\label{prop:G_a-limit}
The operator $G_a(\rho):\calH_p\to \calH_r$ extends continuously
(from $\Cx(X_a)_\alpha^+\cap\Cx(X_a)_\alpha^\circ$) to
$\Cx(X_a)_\alpha^+$ for $p>1/2$, $r<-1/2$, and it satisfies
$P_a(\rho)G_a(\rho)=\Id$, $G_a(\rho)P_a(\rho)=\Id$,
on $\calH_p$, $p>1/2$. The limit $G_a(\nu,\rho_\perp,z\pm i0)$ satisfies
\begin{equation}\begin{split}\label{eq:G_a-real}
&e^{i\rho\cdot w_a}G_a(\nu,\rho_\perp,z\pm i0)e^{-i\rho\cdot w_a}\\
&=R_a(\rho^2+i0)+\sum_{\alpha'}\Fr^{-1}_{X_a}H(\nu\cdot(\xi_a-\re\rho))
[(|\xi_a|^2-(\lambda-\ep_{\alpha'}-i0))^{-1}\\
&\qquad\qquad\qquad\qquad\qquad
-(|\xi_a|^2-(\lambda-\ep_{\alpha'}+i0))^{-1}]\Fr_{X_a}E_{a,\alpha'}.
\end{split}\end{equation}
\end{prop}

The last term in \eqref{eq:G_a-real} can be written in terms of the
Poisson operators in the bound states of $H^a$, using the following
lemma.

\begin{lemma}
Suppose that $\lambda\nin\Lambda'_a$. Then
\begin{equation}\begin{split}\label{eq:G-Poisson}
&\sum_{\alpha'}\Fr^{-1}_{X_a}H(\nu\cdot(\xi_a-\re\rho))\\
&\qquad\qquad\qquad[(|\xi_a|^2-(\lambda-\ep_{\alpha'}-i0))^{-1}
-(|\xi_a|^2-(\lambda-\ep_{\alpha'}+i0))^{-1}]\Fr_{X_a}E_{a,\alpha'}\\
&\qquad\qquad=\sum_{\alpha'}\frac{-i}{2\sqrt{\lambda-\ep_{\alpha'}}}
\calPt_{\alpha'-}(\lambda)H(\nu\cdot(.-\rho))
\calPt_{\alpha'-}(\lambda)^*.
\end{split}\end{equation}
\end{lemma}

\begin{proof}
On $\Range(\Id-E_a)$ both sides vanish.
Thus, it suffices to consider
the equation on $\Range(E_{a,\ep})$, $\ep\in\pspec(H^a)$, or
simply on $\psi_{\alpha'}$ where $\alpha'$ is a bound state of $H^a$ of
energy $\ep=\ep_{\alpha'}$.
Explicitly, on the $X_a$-Fourier transform side, on this space the equation
follows from
\begin{equation*}
(\rho-(\lambda-\ep-i0))^{-1}
-(\rho-(\lambda-\ep+i0))^{-1}
=-2\pi i \delta_{\rho-(\lambda-\ep)}.
\end{equation*}
Indeed, with $\rho=|\xi_a|^2$, hence
$d\xi_a=\frac{1}{2}\rho^{\frac{\dim X_a-2}{2}}\,d\rho\,d\omega_a$,
$u=v\otimes\psi_{\alpha'}$, $\hat v=\Fr_{X_a} v$, this gives
\begin{equation*}\begin{split}
&(2\pi)^{-\dim X_a}\int_{X_a}\left((|\xi_a|^2-(\lambda-\ep-i0))^{-1}
-(|\xi_a|^2-(\lambda-\ep+i0))^{-1}\right)\hat v(w_a)\,dw_a\\
&\qquad=-(2\pi)^{-\dim X_a}\pi i \lambda^{\frac{\dim X_a-2}{2}}
\int_{\sphere_a(\sqrt{\lambda-\ep_\alpha})} e^{iw_a\cdot \omega_a}
\hat v(\omega_a)
\,d\omega_a.
\end{split}\end{equation*}
Comparing with the definition of the Poisson operators, namely that
the Fourier transform, followed by restriction to
$\sphere_a(\sqrt{\lambda-\ep_\alpha})$, is essentially given by
$\calPt_{\alpha'-}(\lambda)^*$ in view of \eqref{eq:P*-F}, with a
similar relation connecting the inverse Fourier transform and
$\calPt_{\alpha'-}(\lambda)$, \eqref{eq:G-Poisson} follows.
\end{proof}

The combination of the preceeding two results yields:

\begin{cor}
Suppose that $\rho\cdot\rho+\ep_\alpha=\lambda
\in(\ep_\alpha,\evth)\setminus\Lambda_a'$,
$\rho_\perp^2+\ep_\alpha\nin\Lambda_a'$. Then
\begin{equation*}\begin{split}
e^{i\rho\cdot w_a}G_a(\rho)e^{-i\rho\cdot w_a}
=&R_a(\lambda+i0)\\
&-\sum_{\alpha'}\frac{i}{2\sqrt{\lambda-\ep_{\alpha'}}}
\calPt_{\alpha'-}(\lambda)
H(\nu\cdot (.-\re\rho))\calPt_{\alpha'-}(\lambda)^*,
\end{split}\end{equation*}
where $\calPt_{\alpha'-}(\lambda)$ is considered as an operator on the sphere
of radius $\lambda-\ep_{\alpha'}$.
\end{cor}

Proposition~\ref{prop:G_a-limit} can be used to show that $G(\rho)$ itself
has a limit when $\rho$ becomes real, provided that
$\lambda=\rho^2+\ep_\alpha<\evth$.

\begin{thm}\label{thm:G-rho-real}
Suppose that $|\rho_\perp|^2+\ep_\alpha\in(\ep_\alpha,\evth)\setminus\Lambda'_a$,
$p>1/2$, $r<-1/2$, $\mu>\max(p,1)$.
There exists $\delta>0$ with the following property.
Suppose that
for all $b\neq a$, $\sup |\langle
w^b\rangle^\mu V_b|<\delta$.
Then
the operator $G(\rho):\calH_p\to \calH_r$ extends continuously to
\begin{equation}\label{eq:limit-region}
\{z:\ \im z>0\}\cup\{z\geq 0:
\ z^2+\rho_\perp^2+\ep_\alpha\in(\ep_\alpha,\evth)\setminus\Lambda'\},
\end{equation}
and it satisfies
$P(\rho)G(\rho)=\Id$, $G(\rho)P(\rho)=\Id$,
on $\calH_p$. 
\end{thm}

\begin{proof}
We only need to show that $(\Id+I_a G_a(\rho))^{-1}$ extends
to real $\rho$ as stated, as a bounded operator on $\calH_p$,
$p>1/2$. But this follows as in the remarks preceeding
Theorem~\ref{thm:G(rho)-exist)}.
\end{proof}

\begin{cor}\label{cor:u_rho-limit}
Suppose that $|\rho_\perp|^2+\ep_\alpha\in(\ep_\alpha,\evth)
\setminus\Lambda'_a$, $\mu>(\dim X_a+1)/2$,
$V_b$ as in Theorem~\ref{thm:G-rho-real}. Then $u_\rho$ extends
continuously to \eqref{eq:limit-region}, with $u_\rho-u^0_\rho\in \calH_r$
for all $r<0$,
and $u_\rho$ is analytic in $z$ in $\im z>0$.
\end{cor}

\section{The connection between the S-matrix and the exponential solutions}

We introduce the analogue of the pairing \eqref{eq:pair-8}
describing the S-matrix via
\begin{equation}\label{eq:pair-16}
G_{\alpha\alpha'}(\rho,\overline{\rho}+\zeta)
=\int_{\Rn} I_a u_{\rho}\overline{u^0_{\alpha',\overline{\rho}+\zeta}}
=\int_{\Rn} I_a u_\rho \overline{\psi_{\alpha'}(w^a)}
e^{-i(\rho+\zeta)\cdot w_a}.
\end{equation}
By Corollary~\ref{cor:u_rho-limit}, if
$\langle w^b\rangle^\mu V_b\in L^\infty(X^b)$ for all $b$
and for some $\mu>\dim X_a$, then the integral in \eqref{eq:pair-16}
converges for all $\rho$ for which $u_\rho$ exists, and for
$\zeta\in\Rn$, since then the real parts of the
exponentials cancel, and
$I_a \psi_\alpha\overline{\psi_{\alpha'}}\in L^1(X_0)$,
and the same holds for $I_a \overline{\psi_{\alpha'}} G(\rho)(I_a\psi_\alpha)$.
Other properties of \eqref{eq:pair-16} follow immediately from
Corollary~\ref{cor:u_rho-limit}.

\begin{prop}
Suppose that $|\rho_\perp|^2+\ep_\alpha\in(\ep_\alpha,\evth)
\setminus\Lambda'_a$, $\mu>\dim X_a$,
and $V_b$ as in
Theorem~\ref{thm:G-rho-real}. Then
$G_{\alpha\alpha'}$ is an
analytic function of $z$ in $\Cx\setminus\Real$, and extends to be
continuous on \eqref{eq:limit-region}.
In addition,
\begin{equation}\label{eq:G-z-infty}
\lim_{|z|\to\infty} G_{\alpha\alpha'}(\rho,\overline{\rho}+\zeta)
=\int_{\Rn} I_a \psi_\alpha \overline{\psi_{\alpha'}} e^{-i\zeta\cdot w_a},
\end{equation}
provided that $|z|\to\infty$ in $|\im z|>C|\re z|$, $C>0$.
\end{prop}

\begin{proof}
The first two statements are direct consequences of
Corollary~\ref{cor:u_rho-limit}. By
Corollary~\ref{cor:u_rho-exists},
$G(\rho)(I_a\psi_\alpha)\to 0$ as $|z|\to\infty$ in
$\calH_r$ for $r<\mu-1-\frac{\dim X_a}{2}<0$. On the other hand,
\begin{equation*}
\overline{\psi_{\alpha'}}I_a e^{-i\zeta\cdot w_a}\in L^2_{s}(X_0)\subset
\calH_s,\ s>0,\ s<\mu-\frac{\dim X_a}{2},
\end{equation*}
hence in $\calH_{-r}$ provided that $r>\frac{\dim X_a}{2}-\mu$. We can take
$r=\mu-1-\frac{\dim X_a}{2}-\ep$, $\ep>0$ sufficiently small. Thus we conclude
that
\begin{equation*}
\lim_{|z|\to\infty} \int_{\Rn} (G(\rho)(I_a\psi_\alpha))
(I_a\overline{\psi_{\alpha'}}e^{-i\zeta\cdot w_a}) =0,
\end{equation*}
hence \eqref{eq:G-z-infty} follows.
\end{proof}

For fixed $\rho$ real and $\zeta$ satisfying
\begin{equation*}
\lambda=\rho^2+\ep_\alpha=(\rho+\zeta)^2+\ep_{\alpha''},
\end{equation*}
i.e.\ the equality of incoming and outgoing energies, we can
relate $G_{\alpha\alpha''}(\rho,\rho+\zeta)$, $\rho$ real, to the
S-matrices as follows.
Under our assumptions,
\begin{equation}\label{eq:R-G-8}
e^{i\rho\cdot w_a}(\Id+G_a(\rho)I_a)e^{-i\rho\cdot w_a}u_\rho
=u^0_\rho=(\Id+R_a(\rho^2+i0)I_a) U_\rho.
\end{equation}
Applying $(\Id+R_a(\rho^2+i0)I_a)^{-1}$ to both sides of
\eqref{eq:R-G-8}, we deduce that
\begin{equation}\begin{split}\label{eq:R-G-16}
U_\rho&=u_\rho-\sum_{\alpha'}\frac{i}{2\sqrt{\lambda-\ep_{\alpha'}}}
(\Id+R_a(\rho^2+i0)I_a)^{-1}\\
&\qquad\qquad\qquad\qquad\calPt_{\alpha'-}(\lambda)H(\nu\cdot (.-\re\rho))
\calPt_{\alpha'-}(\lambda)^*I_a u_\rho\\
&=u_\rho-\sum_{\alpha'}\frac{i}{2\sqrt{\lambda-\ep_{\alpha'}}}
\calP_{\alpha'-}(\lambda)H(\nu\cdot (.-\re\rho))
\calPt_{\alpha'-}(\lambda)^*I_a u_\rho.
\end{split}\end{equation}
Integrating against
$I_a e^{-i(\rho+\zeta)\cdot w_a}\overline{\psi_{\alpha''}(w^a)}$
yields
\begin{equation}\begin{split}
\calS^\sharp_{\alpha\alpha''+}(\lambda,\rho,\rho+\zeta)
=&G_{\alpha\alpha''}(\rho,\rho+\zeta)\\
&-\sum_{\alpha'}
\frac{i}{2\sqrt{\lambda-\ep_{\alpha'}}}
\int_{\sphere_a(\sqrt{\lambda-\ep_{\alpha'}})}
\calS^\sharp_{\alpha'\alpha''+}(\lambda,\rho',\rho+\zeta)\\
&\qquad\qquad\qquad\qquad H(\nu\cdot (\rho'-\re\rho))
G_{\alpha\alpha'}(\rho,\rho')\,d\rho'.
\end{split}\end{equation}
This is an integral equation for $G_{\alpha\alpha''}$ in terms of
$\calS^\sharp_{\alpha\alpha''+}(\lambda)$,
$\calS^\sharp_{\alpha'\alpha''+}(\lambda)$,
$\lambda=\rho^2+\ep_\alpha$.

\begin{prop}\label{prop:G-S}
Suppose that $\rho^2+\ep_\alpha=\lambda
=(\rho+\zeta)^2+\ep_{\alpha''}\in(-\infty,0)\setminus\Lambda'$,
$|\rho_\perp|^2+\ep_\alpha\in(\ep_\alpha,\evth)
\setminus\Lambda'_a$, $\mu>\dim X_a$. There exists
$\delta>0$ with the following property.

Suppose that for all $b\neq a$, $\sup |\langle w^b\rangle^\mu V_b|<\delta$.
Then the pairings $G_{\alpha\alpha''}(\rho,\rho+\zeta)$
are determined by the operators
$\calS^\sharp_{\alpha'\alpha''+}(\lambda)$ given for all
$\alpha'$ and $\alpha''$.
\end{prop}

\begin{proof}
We first discuss the case when the only bound state of
$H^a$ is $\alpha$, or more generally if $\ep_\alpha'>\lambda$ for
$\alpha'\neq\alpha$, i.e.\ the total energy $\lambda$ is just above the
ground state energy. Then we get
\begin{equation}\begin{split}\label{eq:S-G-aa}
\calS^\sharp_{\alpha\alpha+}(\lambda,\rho,\rho+\zeta)
=&G_{\alpha\alpha}(\rho,\rho+\zeta)\\
&-\frac{i}{2\sqrt{\lambda-\ep_{\alpha}}}
\int_{\sphere(\sqrt{\lambda-\ep_{\alpha}})}\calS^\sharp_{\alpha\alpha+}(\lambda,
\rho',\rho+\zeta)H(\nu\cdot (\rho'-\rho))\\
&\hspace*{7cm} G_{\alpha\alpha}(\rho,\rho')\,d\rho'.
\end{split}\end{equation}
Now fix $\rho$, i.e.\ more precisely fix $\nu$, $z$, $\rho_\perp$, and
consider this as an integral equation for the function
$G_{\alpha\alpha}(\rho,.)$. Then this has the form
\begin{equation}
(\Id-T)G_{\alpha\alpha}(\rho,.)=f,
\end{equation}
with $f\in\Cinf(\sphere_a(\sqrt{\lambda-\ep_\alpha}))$, $T:L^2\to C^0$
with bounded kernel and small norm as a map $L^2\to L^2$ since
$\sup_{\omega,\omega'}
|\calS^\sharp_{\alpha\alpha+}(\lambda,\omega,\omega'))|$ is small
by \eqref{eq:pair-8} (as $I_a$ is small). Hence, $\Id-T$
is invertible, proving the proposition in this case.

In complete generality, we consider the vector $\Phi$
whose $\alpha'$ entry is $\Phi_{\alpha'}(.)=G_{\alpha\alpha'}(\rho,.)$.
Then we obtain a system of equations of the form 
\begin{equation}
(\Id-T)\Phi=f,
\end{equation}
as above. Again, $T$ has small norm, so $\Id-T$ is invertible, proving
the proposition.
\end{proof}

Novikov \cite{Novikov:Determination} noticed
that in two-body scattering,
the near-forward values of
$G(\rho,\rho+\zeta)$, i.e.\ the values when the angles between
$\rho$ and $\nu$, resp.\ between $\rho+\zeta$ and $\nu$ are small,
is determined by $\calS^\sharp_{\alpha'\alpha''+}(\rho,\rho')$
where the angle between $\rho$ and $\nu$, resp.\ and $\rho'$ and $\nu$
is small. His observation
also applies in the present setting.

This can be understood via linear algebra. Thus, we decompose
\begin{equation}\begin{split}\label{eq:forward-split}
&V=L^2(\sphere_a(\sqrt{\lambda-\ep}))=V_1\oplus V_2,\\
&V_1=L^2(\{\rho'\in \sphere_a(\sqrt{\lambda-\ep}):\ \rho'\cdot\nu\geq\rho\cdot\nu\}),\\
&V_2=L^2(\{\rho'\in \sphere_a(\sqrt{\lambda-\ep}):\ \rho'\cdot\nu\leq\rho\cdot\nu\}),
\end{split}\end{equation}
writing the two orthogonal projections as $\pi_1$ and $\pi_2$.
Now $T:V\to V$ vanishes on $V_2$, so $T=T\pi_1$, while its
restriction to $V_1$, via $\pi_1$, is exactly $T_1$: $T_1=\pi_1 T$.
Now linear algebra shows that we only need to know $T_1$ and $f_1=\pi_1 f$ to
find $\pi_1(\Id-T)^{-1}f$, namely
\begin{equation}
\pi_1(\Id-T)^{-1}f=(\Id_{V_1}-T_1)^{-1}\pi_1 f.
\end{equation}
This proves the following proposition.

\begin{prop}\label{prop:G-S-2}
Let $\rho$, $\zeta$, $\lambda$, $\mu$, $\delta$ be as in
Proposition~\ref{prop:G-S}, and suppose that
for all $b\neq a$, $\sup |\langle w^b\rangle^\mu V_b|<\delta$.
Then the pairings $G_{\alpha\alpha''}(\rho,\rho+\zeta)$ for $(\rho+\zeta)\cdot
\nu\geq\rho\cdot\nu$ are determined by the S-matrices
$\calS^\sharp_{\alpha'\alpha''+}(\lambda,\omega,\omega')$ given for all
$\alpha'$ and $\alpha''$, evaluated in $\omega\cdot\nu\geq z=\rho\cdot\nu$,
$\omega'\cdot\nu\geq z=\rho\cdot\nu$.
\end{prop}

Note that if $\ep_\alpha$ is the bottom of the spectrum of $H^a$,
$\rho_\perp$ is sufficiently small, namely $|\rho_\perp|^2<\ep'-\ep_\alpha$,
$\ep'$ as in Theorem~\ref{thm:main2},
$\omega'\in\sphere_a(\sqrt{\lambda-\ep_{\alpha'}})$ then
$\omega'\cdot\nu\leq \sqrt{\lambda-\ep_{\alpha'}}$ while $z^2=\lambda
-\rho_\perp^2-\ep_\alpha>\lambda-\ep'$ shows that $\rho\cdot\nu=z
>\sqrt{\lambda-\ep'}$, hence $\omega'\cdot\nu\geq z$ never holds. Thus,
under the conditions of Theorem~\ref{thm:main2}, only
$\calS^\sharp_{\alpha\alpha+}$ is needed to determine $G_{\alpha\alpha}(\rho,
\rho+\zeta)$ for $|\rho_\perp|<\sqrt{\ep'-\ep_\alpha}$, $\zeta\cdot\nu\geq 0$.

\section{Inverse results: proof of Theorems~\ref{thm:main} and \ref{thm:main2}}
In this section we prove Theorems~\ref{thm:main} and \ref{thm:main2}.

Fix a non-empty interval open $I\subset(\ep_\alpha,0)$, and let
$R=2\sqrt{\sup I-\ep_\alpha}$. Let $\zeta\in X_a$ satisfy $|\zeta|<R$, and
$\frac{1}{4}|\zeta|^2+\ep_\alpha\nin\Lambda_a'$. The last condition
excludes a discrete set of values of $|\zeta|^2$ in $[0,R)$.
Note that $\frac{1}{4}|\zeta|^2+\ep_\alpha<0$.

We now choose $\evpr<0$ in \eqref{eq:ep_1-def} so that
$\frac{1}{4}|\zeta|^2+\ep_\alpha<\evpr$, and
\begin{equation}\label{eq:I-energy}
I\cap(\frac{1}{4}|\zeta|^2+\ep_\alpha,\evpr)\neq\emptyset.
\end{equation}
Let $\rho_\perp=-\zeta/2$, so $\rho_\perp^2+\ep_\alpha\in(\ep_\alpha,0)
\setminus\Lambda'_a$.
Let $\nu\in X_a$ be such that $|\nu|^2=1$, $\nu\cdot \zeta=0$.
Then for $\rho=z\nu+\rho_\perp$, $z\in\Cx$,
$(\rho+\zeta)^2-\rho^2=2\rho\cdot\zeta+\zeta^2=0$.
Having fixed $\zeta$, $\rho_\perp$, $\nu$, consider the
energy $\rho^2+\ep_\alpha$ when $z$ becomes real.
Since
$\rho^2=z^2+|\rho_\perp|^2$, as $z$ varies over
$[0,\sqrt{\evth-\rho_\perp^2-\ep_\alpha})$, the
energy varies over $[|\rho_\perp|^2+\ep_\alpha,\evth)$. This
intersects the given interval $I$ by \eqref{eq:I-energy} as $\evpr<\evth$.

Let $z_0>0$ satisfy $\lambda=\rho_\perp^2+z_0^2+\ep_\alpha
\in(\ep_\alpha,\evth)\setminus\Lambda'$.
By Proposition~\ref{prop:G-S}
the limit of the pairing $G_{\alpha\alpha}(\rho,\rho+\zeta)$
as $z\to z_0$ (in $\im z>0$)
is determined by the scattering matrices
$\calS^\sharp_{\alpha'\alpha''}(\rho_\perp^2+z_0^2+\ep_\alpha)$. Thus,
there exists a non-empty open interval of these values $z_0$ at which
$G_{\alpha\alpha}(\rho,\rho+\zeta)$ is determined
by $\calS^\sharp_{\alpha'\alpha''}(\lambda)$, $\lambda\in I$.

Indeed, as explained in and after the statement of
Proposition~\ref{prop:G-S-2}, under the assumptions of
Theorem~\ref{thm:main2}, with $|\zeta|<2\sqrt{\min(\sup I,\ep')-\ep_\alpha}$,
$\ep'$ denoting
the next eigenvalue of $H^a$ or $0$, as in the statement of the theorem,
one only needs to know $\calS^\sharp_{\alpha\alpha}(|\rho|^2+\ep_\alpha)$, and
in either case the knowledge of the S-matrices in appropriate near-forward
regions suffices due the remarks surrounding \eqref{eq:forward-split}.

Since the limit on any open interval in the boundary of its domain
determines an analytic function, we deduce that knowing the S-matrix
$\calS^\sharp_{\alpha'\alpha''}$ in the interval $I$
determines $G_{\alpha\alpha}(\rho,\zeta)$
for all $\zeta$ with $|\zeta|<R$.

Now let $z\to\infty$ through imaginary $z$.
By \eqref{eq:G-z-infty}, $G_{\alpha\alpha}(\rho,\rho+\zeta)$ converges to
\begin{equation}\label{eq:pair-99}
\int_{\Rn} I_a |\psi_{\alpha}|^2 e^{-i\zeta\cdot w_a}
\,dw=\int_{X_a}e^{-i\zeta\cdot w_a} \left(\int_{X^a}
I_a|\psi_{\alpha}|^2\,dw^a\right)\,dw_a.
\end{equation}
But this is the Fourier transform of $\int_{X^a}
I_a|\psi_{\alpha}|^2\,dw^a$ in $X_a$, evaluated at $\zeta$.
Hence $S_{\alpha'\alpha''}(\lambda)$, $\lambda\in I$, determines
the Fourier transform in $w_a$ of the `effective interaction'
\begin{equation}\label{eq:pair-101}
\int_{X^a}
I_a|\psi_{\alpha}|^2\,dw^a
\end{equation}
in a ball of radius $2\sqrt{\sup I-\ep_\alpha}$, except possibly on the spheres
$\frac{1}{4}|\zeta|^2+\ep_\alpha\in\Lambda'_a$. However, $\int_{X^a}
I_a|\psi_{\alpha}|^2\,dw^a\in L^1(X_a)$, hence the Fourier transform is
continuous, hence it is determined on the whole ball.
This completes the proof of Theorems~\ref{thm:main} and
\ref{thm:main2}.

\bibliographystyle{plain}
\bibliography{sm}

\end{document}